\documentclass[12pt]{amsart}

\usepackage{amsmath, amssymb, amsfonts, amsthm, amscd}
\usepackage{manfnt}
\usepackage{mathtools}
\usepackage{fullpage}
\usepackage{url}
\usepackage[all]{xy}
   \SelectTips{cm}{10}
\usepackage{txfonts}
\usepackage{enumitem}
\DeclareFontFamily{U}{wncy}{}
    \DeclareFontShape{U}{wncy}{m}{n}{<->wncyr10}{}
    \DeclareSymbolFont{mcy}{U}{wncy}{m}{n}
    \DeclareMathSymbol{\Sh}{\mathord}{mcy}{"58}

\newtheorem{thm}{Theorem}[section]
\newtheorem{lemma}[thm]{Lemma}

\newtheorem{cor}[thm]{Corollary}

\newtheorem{prop-conj}[thm]{Proposition-Conjecture}

\theoremstyle{definition}
\newtheorem{defn}[thm]{Definition}

\theoremstyle{remark}
\newtheorem{rmk}[thm]{Remark}

\theoremstyle{remark}

\theoremstyle{remark}

\theoremstyle{remark}

\newcommand{\Q}{\mathbb{Q}}

\newcommand{\Z}{\mathbb{Z}}
\newcommand{\CC}{\mathbb{C}}

\newcommand{\Qlb}{\overline{\mathbb{Q}}_\ell}
\newcommand{\af}{\mathbf{A}_F}

\DeclareMathOperator{\Hom}{Hom}

\DeclareMathOperator{\Aut}{Aut}

\DeclareMathOperator{\im}{im}

\DeclareMathOperator{\Fil}{Fil}
\DeclareMathOperator{\gr}{gr}

\DeclareMathOperator{\Rep}{Rep}

\DeclareMathOperator{\id}{id}

\DeclareMathOperator{\Gr}{Gr}

\DeclareMathOperator{\Vect}{Vect}

\newcommand{\gal}[1]{\Gamma_{#1}} 
\newcommand{\Gal}{\mathrm{Gal}} 

\newcommand{\into}{\hookrightarrow}
\newcommand{\onto}{\twoheadrightarrow}

\newcommand{\mc}{\mathcal}
\newcommand{\mf}{\mathfrak}
\newcommand{\mr}{\mathrm}
\newcommand{\mbf}{\mathbf}

\newcommand{\tZ}{\widetilde{Z}}

\newcommand{\tH}{\widetilde{H}}

\newcommand{\oE}{\overline{E}}
\newcommand{\wt}{\widetilde}

\begin{document}
\title{Generalized Kuga-Satake theory and good reduction properties of Galois representations}
\author{Stefan Patrikis}
\email{patrikis@math.utah.edu}
\address{Department of Mathematics\\ The University of Utah\\ Salt Lake City, UT}
\date{April 2016}
\maketitle
\begin{abstract}
In previous work we described when a single geometric representation $\Gamma_F \to H(\overline{\mathbb{Q}}_{\ell})$ of the Galois group of a number field $F$ lifts through a central torus quotient $\widetilde{H} \to H$ to a geometric representation. In this paper we prove a much sharper result for systems of $\ell$-adic representations, such as the $\ell$-adic realizations of a motive over $F$, having common ``good reduction" properties. Namely, such systems admit geometric lifts with good reduction outside a common finite set of primes. The method yields new proofs of theorems of Tate (the original result on lifting projective representations over number fields) and Wintenberger (an analogue of our main result in the case of a central isogeny $\widetilde{H} \to H$). 
\end{abstract}
\section{Introduction}
Let $F$ be a number field, and let $\gal{F}= \Gal(\overline{F}/F)$ be its absolute Galois group with respect to a fixed algebraic closure $\overline{F}$. A fundamental theorem of Tate (\cite[\S 6]{serre:DSsurvey}) asserts that $H^2(\gal{F}, \Q/\Z)$ vanishes; as a result, all (continuous, $\ell$-adic) projective representations of $\gal{F}$ lift to genuine representations, and more generally, whenever $\tH \to H$ is a surjection of linear algebraic groups over $\Qlb$ with kernel equal to a central torus in $\tH$, all representations $\rho_{\ell} \colon \gal{F} \to H(\Qlb)$ lift to $\tH(\Qlb)$.

The $\ell$-adic representations of greatest interest in number theory are those with conjectural connections to the theories of motives and automorphic forms; if the monodromy group of $\rho_{\ell}$ is semi-simple, then it is expected--by conjectures of Fontaine-Mazur, Tate, Grothendieck-Serre, and Langlands--that the $\rho_{\ell}$ arising from pure motives or automorphic forms are precisely those that are geometric in the sense of Fontaine-Mazur, i.e. unramified outside a finite set of places of $F$, and de Rham at all places dividing $\ell$. The paper \cite{stp:variationsarxiv} established a variant of Tate's lifting theorem for such geometric Galois representations. There are subtleties when $F$ has real embeddings, but at least for totally imaginary $F$, any geometric $\rho_{\ell} \colon \gal{F} \to H(\Qlb)$ satisfying a natural ``Hodge symmetry" requirement will admit a geometric lift $\tilde{\rho}_{\ell} \colon \gal{F} \to \tH(\Qlb)$ (see \cite[Theorem 3.2.10]{stp:variationsarxiv}). This geometric lifting theorem leads to a precise expectation for the corresponding lifting problem for motivic Galois representations. Namely, if $\mc{G}_{F, E}$ denotes the motivic Galois group for pure motives over $F$ with coefficients in a number field $E$--we will make this set-up precise in \S \ref{hodgesection}, but for now the reader may take homological motives under the Standard Conjectures--and if $\tH \to H$ is now a surjection of groups over $E$ with central torus kernel, then we conjecture (\cite[Conjecture 4.3.1]{stp:variationsarxiv}) that any motivic Galois representation $\rho \colon \mc{G}_{F, E} \to H$ lifts to $\tH$, at least after some finite extension of coefficients:
\[
\xymatrix{
& \tH_{\overline{E}} \ar[d] \\
\mc{G}_{F, \overline{E}} \ar[r]_{\rho \otimes_E \overline{E}} \ar@{-->}[ur]^{\tilde{\rho}} & H_{\overline{E}}.
}
\]

There is essentially one classical example (with several variants) of this conjecture, a well-known construction of Kuga and Satake (\cite{kuga-satake}), which associates to a complex, for our purposes projective, K3 surface $X$ a complex abelian variety $KS(X)$, related by an inclusion of Hodge-structures $H^2(X, \Q) \subset H^1(KS(X), \Q)^{\otimes 2}$. In the motivic Galois language, finding $KS(X)$ amounts (when $F=\CC$) to finding a lift $\tilde{\rho}$ of the representation $\rho_X \colon \mc{G}_{\CC, \Q} \to H=\mr{SO}(H^2(X)(1))$, through the surjection $\tH= \mr{GSpin}(H^2(X)(1)) \to H$. Progress on the general conjecture, when the motives in question do not lie in the Tannakian sub-category of motives generated by abelian varieties, seems to require entirely new ideas.\footnote{For some, admittedly limited, examples, see \cite{stp:KSRLS1:middleconv} and \cite{stp:KSRLS2:heckeeigen}.} 

The aim of this paper is to establish a Galois-theoretic result which is necessary for this conjecture to hold, but considerably more delicate than the basic geometric lifting theorem of \cite[Theorem 3.2.10]{stp:variationsarxiv}. Namely, any motive $M$ over $F$ has good reduction outside a finite set of primes: for any choice of variety $X$ in whose cohomology $M$ appears, $X$ spreads out as a smooth projective scheme over $\mc{O}_F[1/N]$ for some integer $N$. In particular, by the base-change theorems of \'{e}tale cohomology (\cite{sga4h}) and the crystalline $p$-adic comparison isomorphism (\cite{faltings:crystalline}), for any motivic Galois representation $\rho \colon \mc{G}_{F, E} \to H$, the $\lambda$-adic realizations $\rho_{\lambda} \colon \gal{F} \to H(E_{\lambda})$ have \textit{good reduction} outside a finite set of primes $S$, in the sense (also see Definition \ref{ramcomp}) that each $\rho_{\lambda}$ factors through $\gal{F, S \cup S_{\lambda}}$ and is crystalline at all places of $S_{\lambda} \setminus (S_{\lambda} \cap S)$: here $S_{\lambda}$ denotes the primes of $F$ with the same residue characteristic as $\lambda$, and $\gal{F, S \cup S_{\lambda}}$ is the Galois group of the maximal extension of $F$ inside $\overline{F}$ that is unramified outside of $S \cup S_{\lambda}$. Certainly a necessary condition for the generalized Kuga-Satake conjecture to hold is that the realizations $\{\rho_{\lambda}\}_{\lambda}$ of $\rho$ should lift to geometric representations $\tilde{\rho}_{\lambda} \colon \gal{F, P \cup S_{\lambda}} \to \tH(\oE_{\lambda})$ that likewise have good reduction outside a common finite set of places $P$. This is what we will show, as a consequence of a more general result, whose precise statement we postpone until Theorem \ref{mainintro} (and for the multiplicative-type case, Corollaries \ref{wint} and \ref{multtypecor}). Here is the application in the motivic setting:
\begin{cor}\label{introcor}
Let $F$ be a totally imaginary number field, let $E$ be a number field, and let $\mc{G}_{F, E}$ denote the motivic Galois group, defined by Andr\'{e}'s motivated cycles, of pure motives over $F$ with coefficients in $E$ (see \S \ref{hodgesection}). Let $\tH \to H$ be a surjection of linear algebraic groups over $E$ whose kernel is a central torus in $\tH$, and let $\rho \colon \mc{G}_{F, E} \to H$ be any motivic Galois representation, with associated $\lambda$-adic realizations $\rho_{\lambda} \colon \gal{F, S \cup S_{\lambda}} \to H(\oE_{\lambda})$ for some finite set $S$ of places of $F$. Then there exist a finite, independent of $\lambda$, set $P \supset S$ of places of $F$ and, for all $\lambda$, lifts
\[
\xymatrix{
& \tH(\oE_{\lambda}) \ar[d] \\
\gal{F, P \cup S_{\lambda}} \ar[r]_{\rho_{\lambda}} \ar@{-->}[ur]^{\tilde{\rho}_{\lambda}} & H(\oE_{\lambda})
}
\]
such that each $\tilde{\rho}_{\lambda}$ is de Rham at all places in $S_{\lambda}$, and is moreover crystalline at all places in $S_{\lambda} \setminus (S_{\lambda} \cap P)$.

Now suppose $H' \to H$ is a surjection of linear algebraic groups whose kernel is central but of arbitrary multiplicative type, and let $\rho \colon \mc{G}_{F, E} \to H$ again be a motivic Galois representation. Assume that the labeled Hodge co-characters of $\rho$ (see Definition \ref{hodgeco}) lift to $H'$. Then there exist a finite set $P \supset S$ of places of $F$, a finite extension $F'/F$, and, for all $\lambda$, lifts 
\[
\xymatrix{
& H'(\oE_{\lambda}) \ar[d] \\
\gal{F', P \cup S_{\lambda}} \ar[r]_{\rho_{\lambda}|_{\gal{F'}}} \ar@{-->}[ur]^{\tilde{\rho}_{\lambda}} & H(\oE_{\lambda})
}
\]
such that $\tilde{\rho}_{\lambda}$ is de Rham at all places above $S_{\lambda}$, and is moreover crystalline at all places above $S_{\lambda} \setminus (S_{\lambda} \cap P)$.
\end{cor}
\begin{rmk}\label{real}
All results of this paper, once we take into account the caveat of \cite[\S 2.8]{stp:variationsarxiv} (see too \cite[Proposition 5.5]{stp:parities}), admit straightforward variants when $F$ has real places. Thus for real $F$, the analogue of Corollary \ref{introcor} either holds exactly as written, or after replacing $F$ by any totally imaginary (eg, composite with a quadratic imaginary) extension. We do not want to discuss this at any length here, but we simply remind the reader that the prototypical example in which $F$ is totally real, and Corollary \ref{introcor} fails as stated, is that of the projective motivic Galois representation associated to a \textit{mixed-parity} Hilbert modular form.
\end{rmk}
The isogeny case of this theorem (Corollary \ref{wint}) is a beautiful result of Wintenberger (\cite[Th\'{e}or\`{e}me 2.1.4, Th\'{e}or\`{e}me 2.1.7]{wintenberger:relevement}), and we note here that the finite base-change $F'/F$ may in that case be necessary: consider for instance the projectivization of the Tate module of an elliptic curve over $\Q$, where the $\mr{PGL}_2$-valued Galois representation only lifts to $\mr{SL}_2$ after a finite base-change (over which the $\ell$-adic cyclotomic character acquires a square root). Our problem resembles Wintenberger's in that both lead to a basic difficulty of annihilating cohomological obstruction classes in infinitely many Galois cohomology groups, one for each $\lambda$, but needing to do so in an ``independent-of-$\lambda$" fashion. The arguments themselves, however, are in fact orthogonal to one another: Wintenberger always kills cohomology by making a finite base-change on $F$, whereas that is precisely what we are forbidden from doing if we want the more precise results of Corollary \ref{introcor}. In fact, our theorem in the case of central torus quotients readily implies Wintenberger's main theorem (and its generalization to quotients of multiplicative type; see \S \ref{multtype}); and our arguments also yield a novel proof of Tate's original vanishing theorem (see Corollary \ref{tate}). In fact, Corollary \ref{tate} establishes a more precise form of Tate's theorem: the latter of course shows that the image under the canonical map $H^2(\gal{F, S}, \Z/N) \to H^2(\gal{F}, \Q/\Z)$ is zero, and our refinement quantifies how much additional ramification must be added, and how much the coefficients must be enlarged, in order to annihilate $H^2(\gal{F, S}, \Z/N)$. Our arguments thus achieve, from scratch, a satisfying common generalization of the theorems of Wintenberger and Tate.

In Corollaries \ref{charcomp} and \ref{abcomp}, we give a couple of applications to lifting $\lambda$-adic realizations such that the associated ``similitude characters'' (eg, determinant or Clifford norm) of the lifts form strongly-compatible systems. Note that even in the case of the classical Kuga-Satake construction, this compatibility is only achieved as a consequence of having an arithmetic descent of the (Hodge-theoretically defined) Kuga-Satake abelian variety; such a descent depends on the deformation theory of K3 surfaces and monodromy arguments (due to Deligne and Andr\'{e}: see \cite{deligne:weilK3} and \cite{andre:hyperkaehler}). 

Finally, let us emphasize what we do \textit{not} prove. The realizations $\{\rho_{\lambda}\}_{\lambda}$ of $\rho$ should moreover form a weakly-compatible system of Galois representations in the sense that the conjugacy class of $\rho_{\lambda}(fr_v)$ is defined over $\oE$ and is suitably independent of $\lambda$ (for $v$ outside $S \cup S_{\lambda}$), and in turn one would hope to construct lifts $\tilde{\rho}_{\lambda}$ with the same frobenius compatibility. This problem seems to be out of reach: I know of no way to establish such results using only Galois-theoretic techniques, although indeed they would follow (assuming the Standard Conjectures) from the generalized Kuga-Satake conjecture.\footnote{Alternatively, it is possible to establish results of this nature in settings where $\rho_{\lambda}$ and $\tilde{\rho}_{\lambda}$ are constructed as automorphic Galois representations.} 
\subsection{Notation}\label{notation}
For a number field $F$, we fix an algebraic closure $\overline{F}$ and set $\gal{F}= \Gal(\overline{F}/F)$. For a finite set of places $S$ of $F$, we let $F(S)$ denote the maximal extension of $F$ inside $\overline{F}$ that is unramified outside of $S$, and then we set $\gal{F, S}= \Gal(F(S)/F)$. We denote the ring of $S$-integers in $F$ by $\mc{O}_F[\frac{1}{S}]$. We also set $F_S= \prod_{v \in S} F_v$. If $L$ is a finite extension of $F$ (inside $\overline{F}$), we then abusively continue to write $S$ for the set of all places of $L$ above those in $S$, with corresponding notation $L(S)$, $\gal{L, S}$, etc. 
\section{Hodge symmetry}\label{hodgesection}
In this section we establish a motivic setting in which our general Galois-theoretic results apply; this setting will both serve as motivation for subsequent sections and allow us to deduce Corollary \ref{introcor} from our main Theorem \ref{main} (and Corollary \ref{multtypecor}). Rather than working with (pure) homological motives and assuming the Standard Conjectures, we work with a category of motives that is unconditionally semi-simple and Tannakian--and in which we can prove unconditional results--but that would, under the Standard Conjectures, turn out to be equivalent to the category of homological motives. Namely, let $\mc{M}_{F, E}$ denote Andr\'{e}'s category of motivated motives over $F$ with coefficients in $E$ (see \cite{andre:motivated}). We begin by elaborating on the consequences of Hodge symmetry in $\mc{M}_{F, E}$. Throughout this discussion, it will be convenient to fix embeddings $\tau_0 \colon F \into \oE$ and $\iota_{\infty} \colon \oE \into \CC$. The composite $\iota_{\infty} \tau_0 \colon F \into \CC$ yields a Betti fiber functor $H_{\iota_{\infty}\tau_0} \colon \mc{M}_{F, E} \to \Vect_E$, making $\mc{M}_{F, E}$ into a neutral Tannakian category over $E$. We denote by $\mc{G}= \mc{G}_B(\iota_{\infty} \tau_0)$ the associated Tannakian group (tensor automorphisms of the fiber functor), so that $H_{\iota_{\infty} \tau_0}$ induces an equivalence of tensor categories $\mc{M}_{F, E} \xrightarrow{\sim} \Rep(\mc{G})$. 

We will consider other cohomological realizations on $\mc{M}_{F, E}$, and their comparisons with the Betti fiber functor. Let $H_{dR} \colon \mc{M}_{F, E} \to \Fil_{F \otimes_{\Q} E}$ denote the de Rham realization, taking values in filtered $F \otimes_{\Q} E$-modules; and for each place $\lambda$ of $E$, let $H_{\lambda}$ denote the $\lambda$-adic realization, which takes values in finite $E_{\lambda}$-modules with a continuous action of $\gal{F}$. For all embeddings $\tau \colon F \into \oE$, we obtain an $\oE$-valued fiber functor
\[
\omega_{dR, \tau} \colon M \mapsto \gr^{\bullet} \left( e_{\tau} H_{dR}(M) \right),
\]
where $e_{\tau}$ is the idempotent induced by $\tau \otimes 1 \colon F \otimes_{\Q} E \to \oE$. Let $\mc{G}_{dR}(\tau)= \Aut^{\otimes}(\omega_{dR, \tau})$ be the associated Tannakian group over $\oE$. Of course this fiber functor factors through the category $\Gr_{\oE}$ of graded $\oE$-vector spaces,
\[
\xymatrix{
\mc{M}_{F, \oE} \ar[dr] \ar[rr] & &\Gr_{\oE} \ar[dl]\\
& \Vect_{\oE}, & 
}
\]
so we obtain a corresponding homomorphism $\mu_{\tau} \colon \mathbf{G}_{m, \oE} \to \mc{G}_{dR}(\tau)$. Without specifying $\tau$, we obtain a fiber functor (see \cite[\S 3]{deligne-milne}) $\omega_{dR}= \gr^{\bullet} H_{dR}$ on $\mc{M}_{F, E}$ valued in projective $F \otimes_{\Q} E$-modules. By \cite[Theorem 3.2]{deligne-milne}, the functor $\underline{\Hom}^{\otimes}(H_{\iota_{\infty}\tau_0}, \omega_{dR})$ is a $\mc{G}$-torsor over $F \otimes_{\Q} E$. In particular, for all $\tau \colon F \into \oE$, we can choose a point of $\underline{\Hom}^{\otimes}(H_{\iota_{\infty} \tau_0} \otimes_E \oE, \omega_{dR, \tau})$ to induce a co-character $\mu_{\tau}$ of $\mc{G}_{\oE}$; and the conjugacy class $[\mu_{\tau}]$ of $\mu_{\tau}$ is independent of this choice.
\begin{defn}\label{hodgeco}
For each $\tau \colon F \into \oE$, we call any $\mu_{\tau} \colon \mbf{G}_{m, \oE} \to \mc{G}_{\oE}$ as above a $\tau$-labeled Hodge co-character; it is a representative of the conjugacy class of co-characters $[\mu_{\tau}]$, the latter being canonically independent of any of the above choices of isomorphisms of fiber functors.
\end{defn}

Next we observe:
\begin{lemma}\label{cm}
For all $\sigma \in \Gal(\oE/E)$, $[\mu_{\tau}]= [\mu_{\sigma \tau}]$. In particular, $[\mu_{\tau}]$ only depends on the restriction of $\tau$ to the maximal CM (or totally real) subfield $F_{cm}$ of $F$.
\end{lemma}
\proof
We decompose $F \otimes_{\Q} E = \prod_i E_i$ into a product of fields, writing $p_i$ for the projection onto $E_i$. Any $E$-algebra homomorphism $\tau \colon F \otimes_{\Q} E \to \oE$ factors through $p_{i(\tau)}$ for a unique $i(\tau)$, and then the $\Gal(\oE/E)$-orbit of $\tau$ is precisely those $E$-algebra homomorphisms (i.e., embeddings $F \into \oE$) $\tau' \colon F \otimes_{\Q} E \to \oE$ such that $i(\tau)= i(\tau')$. The first claim follows, since both $\omega_{dR, \tau}$ and $\omega_{dR, \sigma \tau}$ can be factored through $p_{i(\tau)} \circ \omega_{dR}$. The second claim follows from the first, and the fact that all motives arise by scalar extension from motives with coefficients in CM (or totally real) fields: see \cite[Lemma 4.1.22]{stp:variationsarxiv}.
\endproof

Next note that the canonical weight-grading on $\mc{M}_{F, E}$ induces a \textit{central} weight homomorphism
\[
\omega \colon \mbf{G}_{m, E} \to \mc{G},
\]
and likewise for any other choice of fiber functor and Tannakian group (because $\omega$ is central, it is in fact canonically independent of any choice of isomorphism between fiber functors). Hodge symmetry then results from the complex conjugation action on Betti cohomology, interpreted via the Betti-de Rham comparison isomorphism, which is a distinguished $\CC$-point of $\underline{\Hom}^{\otimes}(\omega_{dR, \tau} \otimes_{\oE, \iota_{\infty}} \CC, H_{\iota_{\infty} \tau} \otimes_{E, \iota_{\infty}} \CC)$.
Namely, complex conjugation on complex-analytic spaces induces (see \cite[Lemma 4.1.24]{stp:variationsarxiv}) natural isomorphisms (without restricting to particular graded pieces for the weight and Hodge filtrations, these are isomorphisms of fiber functors over $\CC$)
\begin{equation}\label{Finfty}
\gr^p \left( e_{\tau} H^w_{dR}(M) \right) \otimes_{\oE, \iota_{\infty}} \CC \xrightarrow{\sim} \gr^{w-p} \left( e_{c \tau} H^w_{dR}(M) \otimes_{\oE, \iota_{\infty}} \CC \right),
\end{equation}
where $c \in \Aut(\oE)$ is the choice of complex conjugation for which $\overline{\iota_{\infty}} \tau= \iota_{\infty} c \tau$. We deduce the following relation: 
\begin{lemma}\label{hodgesym}
For any embedding $\tau \colon F \into \oE$, and any choice of complex conjugation $c \in \Aut(\oE)$, the conjugacy classes of co-characters $[\mu_{\tau}]$ and $[\mu_{c \tau}]$ satisfy
\[
[\mu_{\tau}]= \omega \cdot [\mu_{c \tau}^{-1}],
\]
where recall $\omega$ is the weight co-character.
\end{lemma}
\proof
For the choice of complex conjugation specified by $\overline{\iota_{\infty}} \tau= \iota_{\infty} c \tau$, the relation $[\mu_{\tau}]= \omega \cdot [\mu_{c \tau}^{-1}]$ follows, after base-extension $\iota_{\infty} \colon \oE \to \CC$, from Equation (\ref{Finfty}) above; but this relation necessarily descends to $\oE$, since the conjugacy classes of co-characters are defined over any algebraically closed sub-field of $\CC$. It only remains to observe that $[\mu_{c \tau}]$ is independent of the choice of complex conjugation on $\oE$. This follows from the second assertion of Lemma \ref{cm}.
\endproof

The comparison isomorphisms of $p$-adic Hodge theory then imply that the analogue of Lemma \ref{hodgesym} also holds for the associated Hodge-Tate co-characters. For any place $\lambda$ of $E$, fix an algebraic closure $\oE_{\lambda}$. Embeddings $\tau \colon F \into \oE$ and $\iota_{\lambda} \colon \oE \into \oE_{\lambda}$ then induce $\tau_{\iota_{\lambda}} \colon F_v \into \oE_{\lambda}$ for a suitable place $v$ of $F$ of the same residue characteristic $p$ as $\lambda$. 
Meanwhile, the restriction to $\gal{F_v}$ of the $\lambda$-adic realization induces
\begin{equation}\label{HTcochar}
\xymatrix{
\mc{M}_{F, E} \ar[r]^-{H_{\lambda}|_{\gal{F_v}}} & \Rep^{dR}_{E_{\lambda}}(\gal{F_v}) \ar[r]^{D_{dR}} & \Fil_{F_v \otimes_{\Q_p} E_{\lambda}} \ar[r]^{e_{\tau_{\iota_{\lambda}}}} & \Fil_{\oE_{\lambda}} \ar[r]^{\gr} & \Gr_{\oE_{\lambda}} \ar[r] & \Vect_{\oE_{\lambda}},
}
\end{equation} 
where $D_{dR} \colon \Rep_{E_{\lambda}}^{dR}(\gal{F_v}) \to \Fil_{F_v \otimes_{\Q_p} E_{\lambda}}$ denotes Fontaine's functor restricted to the category of de Rham representations. (Here we have invoked Faltings' $p$-adic de Rham comparison isomorphism from \cite{faltings:crystalline}, and the fact--already noted by Andr\'{e} in \cite{andre:motivated}--that it extends to a comparison isomorphism on all of $\mc{M}_{F, E}$; for details of the latter point, see \cite[Lemma 4.1.25]{stp:variationsarxiv}.) Of course, $\Rep^{dR}_{E_{\lambda}}(\gal{F_v})$ also has its standard forgetful fiber functor (let us say $\oE_{\lambda}$-valued), yielding a Tannakian group $\Gamma^{dR}_{v, \lambda}$ for de Rham $\gal{F_v}$-representations over $\oE_{\lambda}$; by choosing an isomorphism between the two $\oE_{\lambda}$-valued fiber functors on $\Rep^{dR}_{E_{\lambda}}(\gal{F_v})$, we obtain a canonical conjugacy class (recall \cite[Theorem 3.2]{deligne-milne}) of ``$\tau_{\iota_{\lambda}}$-labelled Hodge-Tate co-characters'' $[\mu_{\tau_{\iota_{\lambda}}}]$ of $\Gamma^{dR}_{v, \lambda}$. Specializing, this construction defines the labelled Hodge-Tate co-characters of any de Rham Galois representation $\rho \colon \gal{F_v} \to H(E_{\lambda})$, for any affine algebraic group $H$ over $E_{\lambda}$. 

To relate the $\tau_{\iota_{\lambda}}$-labelled Hodge-Tate co-characters in the motivic setting to the Hodge co-characters previously discussed, note that the de Rham comparison isomorphism (\cite{faltings:crystalline}) yields a natural isomorphism of tensor functors $\mc{M}_{F, E} \to \Gr_{\oE_{\lambda}}$,
\[
\gr \left( e_{\tau} \left( H_{dR}(M) \otimes_{E} \oE \right) \otimes_{\oE, \iota_{\lambda}} \oE_{\lambda} \right)\cong \gr \left( e_{\tau_{\iota_{\lambda}}} \left( D_{dR}(H_{\lambda}(M)|_{\gal{F_v}}) \otimes_{E_{\lambda}} \oE_{\lambda} \right) \right).
\]
We deduce:
\begin{cor}
For any embedding $\tau \colon F \into \overline{E}$, and any embedding $\iota_{\lambda} \colon \oE \into \oE_{\lambda}$, there is an equality of conjugacy classes
\[
[\mu_{\tau} \otimes_{\oE, \iota_{\lambda}} \oE_{\lambda}]= [\mu_{\tau_{\iota_{\lambda}}}].
\]
In particular, for all $\lambda$, and for all $E$-embeddings $\iota_{\lambda} \colon \oE \into \oE_{\lambda}$, the conjugacy classes $[\mu_{\tau_{\iota_{\lambda}}}]$ are independent of $(\lambda, \iota_{\lambda})$ when regarded as valued in the common group $\mc{G}_{\oE}$.
\end{cor}
\section{Lifting}\label{lifting}
In this section we prove our main results, all of which will follow from Theorem \ref{mainintro}, stated at the end of this subsection. First we recall the setting and formalize some terminology. For a place $\lambda$ of $E$, let $S_{\lambda}$ denote the set of places of $F$ with the same residue characteristic as $\lambda$. Let $\tH \to H$ be a surjection of linear algebraic groups over $E$ whose kernel is a central torus in $\tH$; eventually we will also consider the case of a central kernel of multiplicative type (eg, a central isogeny), but the discussion in that case will be a straightforward consequence of the torus case. For convenience in formulating our results, we introduce some (non-standard) terminology:
\begin{defn}\label{ramcomp}
A collection $\{\rho_{\lambda} \colon \gal{F} \to H(\oE_{\lambda})\}_{\lambda}$, as $\lambda$ varies over finite places of $E$, of geometric Galois representations is \textit{ramification-compatible} if there exist
\begin{enumerate}
\item a finite set $S$ of places of $F$ such that each $\rho_{\lambda}$ is unramified outside of $S \cup S_{\lambda}$, i.e. factors through
\[
\rho_{\lambda} \colon \gal{F, S \cup S_{\lambda}} \to H(\oE_{\lambda}),
\]
and for $v$ in  $S_{\lambda}$ but not in $S$, $\rho_{\lambda}|_{\gal{F_v}}$ is crystalline;
and
\item a central co-character $\omega \colon \mathbf{G}_{m, E} \to H$ and a collection of conjugacy classes 
\[
\left\{ [\mu_{\tau} \colon \mathbf{G}_{m, \oE} \to H_{\oE}]\right\}_{\tau \colon F \into \oE}
\]
satisfying $[\mu_{\tau}]= \omega \cdot [\mu_{c \tau}^{-1}]$ for any choice of complex conjugation $c \in \Gal(\oE/\Q)$, such that for all $E$-embeddings $\iota_{\lambda} \colon \oE \into \oE_{\lambda}$, inducing via $\tau$ some $\tau_{\iota_{\lambda}} \colon F_v \into \oE_{\lambda}$, the conjugacy class $[\mu_{\tau} \otimes_{\oE, \iota_{\lambda}} \oE_{\lambda}]$ is equal to the conjugacy class of $\tau_{\iota_{\lambda}}$-labelled Hodge-Tate co-characters associated to $\rho_{\lambda}|_{\gal{F_v}}$.
\end{enumerate} 
If a single representation $\rho_{\lambda}$ satisfies the condition in item (1), we will say $\rho_{\lambda}$ \textit{has good reduction outside} $S$.
\end{defn}
\begin{rmk}
\begin{itemize}
\item Note that the $\rho_{\lambda}$ need not be `compatible' in the usual sense (frobenii acting compatibly): if the ``coefficients'' of the $\rho_{\lambda}$ are bounded in a rather strong sense--there exists a common number field over which their frobenius characteristic polynomials are defined--one expects that our collection of $\rho_{\lambda}$ should partition (dividing up the $\lambda$'s) into finitely many compatible systems.\footnote{To see the relevance of bounding the coefficients, the reader may contrast the case of elliptic curves (over $\Q$, say) unramified outside $S$ with that of all weight-two modular forms unramified outside $S$: of the former there are finitely many isogeny classes, since the conductor is bounded, whereas the latter can have level divisible by arbitrarily high powers of the primes in $S$.}
\item The Hodge symmetry requirement of part (2) of Definition \ref{ramcomp} is not the most general constraint that pertains to a compatible system of $\ell$-adic representations. It is a natural condition coming from algebraic geometry, as we saw in Lemma \ref{hodgesym}, but note that there the ambient group was always taken to be the motivic Galois group $G_{\rho}$ of the given motivic Galois representation $\rho$, and not some larger group. But there may be compatible systems (of motivic origin) where the criterion in part (2) of Definition \ref{ramcomp} fails; for example, consider $\rho_{\ell} \colon \gal{\Q} \to \mr{PGL}_3(\Qlb)$ given by the projectivization of $\kappa_{\ell} \oplus 1 \oplus 1$, with $\kappa_{\ell}$ denoting the $\ell$-adic cyclotomic character. The relation $[\mu_{\tau}]= \omega \cdot [\mu_{c \tau}^{-1}]$ implies in this case $[\mu_\tau]= [\mu_{\tau}^{-1}]$, which is false. A way around this is given in \cite[\S 3.2]{stp:variationsarxiv}, where the Hodge symmetry hypothesis is formulated in a way that conjecturally holds for any geometric representation $\rho_{\ell} \colon \gal{F} \to H(\Qlb)$, regardless of its (reductive) algebraic monodromy group. The proof of \cite[Theorem 3.2.10]{stp:variationsarxiv} thus requires a slightly trickier group-theoretic argument than the one we require here. We have opted in this paper to keep the simpler condition (2) above, so as to focus on what is new in the arguments, and because of its obvious centrality from a motivic point of view (in particular, its sufficiency for Corollary \ref{introcor}).
\end{itemize}
\end{rmk}

In the introduction we stated a particularly palatable corollary of the main result of this paper. We now state the main result precisely and take a moment to explain what is gained in the more general phrasing:  
\begin{thm}\label{mainintro}
Let $E$ be a number field, and let $\tH \to H$ be a surjection of linear algebraic groups over $E$ with kernel equal to a central torus in $\tH$. Let $F$ be a totally imaginary number field, and let $S$ be a finite set of places of $F$ containing the archimedean places. Fix a set of co-characters $\{\mu_{\tau}\}_{\tau \colon F \into \oE}$ as in part (2) of Definition \ref{ramcomp}. Then there exists a finite set of places $P \supset S$ such that for any place $\lambda$ of $E$, any embedding $\iota_{\lambda} \colon \oE \into \oE_{\lambda}$, and any geometric representation $\rho_{\lambda} \colon \gal{F, S \cup S_{\lambda}} \to H(\oE_{\lambda})$ such that
\begin{itemize} 
\item $\rho_{\lambda}$ has good reduction outside $S$; and
\item the conjugacy classes of labelled Hodge-Tate co-characters of $\rho_{\lambda}$ are induced via $\iota_{\lambda}$ from $\{\mu_{\tau}\}_{\tau \colon F \into \oE}$,
\end{itemize}
the representation $\rho_{\lambda}$ admits a geometric lift $\tilde{\rho}_{\lambda} \colon \gal{F, P \cup S_{\lambda}} \to \tH(\oE_{\lambda})$ having good reduction outside $P$.

If $\{\rho_{\lambda} \colon \gal{F, S \cup S_{\lambda}} \to H(\oE_{\lambda})\}_{\lambda}$ is a ramification-compatible system, then there exist a finite set of places $P \supset S$ and lifts $\tilde{\rho}_{\lambda} \colon \gal{F, P \cup S_{\lambda}} \to \tH(\oE_{\lambda})$ such that $\{\tilde{\rho}_{\lambda}\}_{\lambda}$ is a ramification-compatible system.
\end{thm}
The proof of this theorem is completed in Theorem \ref{main}. The typical application of this theorem would be to the collection of Galois representations $\{\rho_{\lambda}\}_{\lambda}$ associated to a motivic Galois representation $\rho \colon \mc{G}_{F, E} \to H$, and indeed (the torus quotient case of) Corollary \ref{introcor} of the Introduction follows immediately (using Lemma \ref{hodgesym}). But even admitting a strong finiteness conjecture, that there are finitely many isomorphism classes of such $\rho$, having coefficients in $E$, prescribed Hodge-Tate co-characters, and good reduction outside a fixed finite set $S$, this theorem still says rather more, since even for fixed $\lambda$ it applies to infinitely many distinct $\rho_{\lambda}$ simultaneously (because we have not bounded the coefficients: recall the example of modular forms of weight two whose nebentypus characters have unbounded conductor, even though supported on the fixed finite set $S$ of primes).

For the sake of this additional generality, I will always state results in the form of Theorem \ref{mainintro}, with the most general version followed by the application to ramification-compatible systems (and therefore to $\lambda$-adic realizations of motivic Galois representations).
\subsection{Torus quotients}
In this subsection we prove Theorem \ref{mainintro}. We begin in the next few paragraphs by gathering together all of the ``independent of $\lambda$ and $\rho_{\lambda}$'' data, and the auxiliary constructions we make on top of this data. Fix as in the statement of Theorem \ref{mainintro} a totally imaginary field $F$ (see Remark \ref{real}), a finite set of places $S$ of $F$ containing the infinite places, a number field $E$, and a surjection $\tH \to H$ of linear algebraic groups over $E$ whose kernel is a central torus, which we denote by $C^{\vee}$. (We will without comment also write $C^{\vee}$ for the base-change to various algebraically-closed fields containing $E$.) Next fix an isogeny-complement $H_1$ of $C^{\vee}$ in $\tH$ (for existence of such $H_1$, see \cite[Proposition 5.3, Step 1]{conrad:dualGW}): thus $H_1 \cdot C^{\vee}= \tH$, and $H_1 \cap C^{\vee}$ is finite. For technical reasons, we will later want to include all primes dividing $\#(H_1 \cap C^{\vee})(\oE)$ in the set of bad primes $S$; this will be indicated at the necessary point (see the discussion following Lemma \ref{loctriv}), but it does no harm simply to add these primes to $S$ from now. Consider the quotient map $\nu \colon \tH \to \tH/H_1$; $\tZ^\vee= \tH/H_1$ is a torus, 
and there is an isogeny $C^{\vee} \to \tZ^\vee$, with kernel $C^{\vee} \cap H_1$. Fix a split torus $\tZ$ over $F$ whose dual group (constructed over $\oE$) is isomorphic to $\tZ^\vee \otimes_E \oE$, and fix such an identification (implicit from now on).  

Fix a set of co-characters $\{\mu_{\tau} \colon \mathbf{G}_{m, \oE} \to H_{\oE}\}_{\tau \colon F \into \oE}$, and a central co-character $\omega \colon \mathbf{G}_{m, \oE} \to H_{\oE}$, satisfying the Hodge symmetry requirement of item (2) of Definition \ref{ramcomp}. Denote by $F_{cm}$ the maximal CM subfield of $F$. The condition in Definition \ref{ramcomp} implies that the co-character $\mu_{\tau}$ depends only on the restriction of $\tau$ to $F_{cm}$; we denote this restriction by $\tau_{cm} \colon F_{cm} \into \oE$. We fix a set of representatives $I$ of $\Hom(F_{cm}, \oE)$ modulo complex conjugation, and for each $\sigma \in I$, we fix a lift $\widetilde{\mu_{\sigma}}$ to $\tH$ of $\mu_{\sigma}$, as well as a (central) lift $\wt{\omega}$ of $\omega$. Note that this is possible, because $C^{\vee}$ is a torus. If $\tau \colon F \into \oE$ restricts to a $\sigma \in I$, we then set $\wt{\mu_{\tau}}= \wt{\mu_{\sigma}}$; if not, then $c \tau \colon F \into \oE$ restricts to a $\sigma \in I$, and we then set $\wt{\mu_{\tau}}= \wt{\omega} \wt{\mu_{\sigma}}^{-1}$.

\begin{lemma}\label{hecke}
Fix once and for all an embedding $\iota_{\infty} \colon \oE \into \CC$. There exists an algebraic automorphic representation $\psi$ of $\tZ(\af)$ such that for all $\tau \colon F \into \oE$, inducing $\tau_{\iota_{\infty}} \colon F_v \into \CC$ by composition with $\iota_{\infty}$, the local component $\psi_v \colon F_v^\times \to \CC^\times$ is given by
\[
\psi_v(z)= \tau_{\iota_{\infty}}(z)^{\nu \circ \wt{\mu_{\tau}}} \overline{\tau_{\iota_{\infty}}}(z)^{\nu \circ \wt{\mu_{c \tau}}}.
\]
(Recall that $\nu$ is the quotient $\tH \to \tZ^\vee$.)
\end{lemma}
\proof
We readily reduce to the case $\tZ= \mathbf{G}_m$, where it follows from the description, due to Weil, of the possible archimedean components of algebraic Hecke characters (see \cite{weil:characters}). (This is where Hodge-symmetry is required.)
\endproof
From now on we fix such a $\psi$, and we let $T$ denote the finite set of places of $F$ such that $\psi$ is unramified outside $T$. For any embedding $\iota_{\lambda} \colon \oE \into \oE_{\lambda}$, we can then consider the $\lambda$-adic realization\footnote{To be precise, this depends on $\iota_{\infty}$ and $\iota_{\lambda}$; but $\iota_{\infty}$ is fixed throughout the paper.}
\[
\psi_{\iota_{\lambda}} \colon \gal{F, T \cup S_{\lambda}} \to \tZ^\vee(\oE_{\lambda}).
\]
Each $\psi_{\iota_{\lambda}}$ is a geometric Galois representation, with good reduction outside $T$, and for any $\tau \colon F \into \oE$, inducing $\tau_{\iota_{\lambda}} \colon F_v \into \oE_{\lambda}$, the Hodge-Tate co-character of $\psi_{\iota_{\lambda}}$ associated to $\tau_{\iota_{\lambda}}$ is $\nu \circ \wt{\mu_{\tau}} \otimes_{\oE, \iota_{\lambda}} \oE_{\lambda}$.

Now we consider any geometric representation
\[
 \rho_{\lambda} \colon \gal{F, S \cup S_{\lambda}} \to H(\oE_{\lambda})
\]
having good reduction outside $S$, along with an embedding $\iota_{\lambda} \colon \oE \into \oE_{\lambda}$ such that the Hodge-Tate co-characters of $\rho_{\lambda}$ arise from the collection $\{\mu_{\tau} \colon \mathbf{G}_{m, \oE} \to H_{\oE}\}_{\tau \colon F \into \oE}$ via $\iota_{\lambda}$. Because the kernel of $\tH \to H$ is a central torus, a fundamental theorem of Tate (\cite[\S 6]{serre:DSsurvey}) ensures in this case that $\rho_{\lambda}$, \textit{as a representation of} $\gal{F}$, lifts to $\tH$. 
As we will see, our arguments in fact imply Tate's theorem (Corollary \ref{tate}), so we do not need to assume it in what follows.

We can define an obstruction class $\mc{O}(\rho_{\lambda})$ to lifting $\rho_{\lambda}$ to a continuous representation $\gal{F, S \cup S_{\lambda}} \to H_1(\oE_{\lambda})$ in the usual way: choose a topological (but not group-theoretic) lift $\rho'_{\lambda}$, and then form the 2-cocycle $(g,h) \mapsto \rho'_{\lambda}(gh) \rho'_{\lambda}(h)^{-1}\rho'_{\lambda}(g)^{-1}$, defining 
\[
\mc{O}(\rho_{\lambda}) \in H^2(\gal{F, S \cup S_{\lambda}}, H_1 \cap C^{\vee}).
\]
Here and in what follows, we simply write $H_1 \cap C^\vee$ for the $\oE_{\lambda}$-points of this finite group scheme.
\begin{rmk}\label{thepoint}
Here lies the essential difficulty to be overcome: while Tate's theorem allows us to annihilate the cohomology classes $\mc{O}(\rho_{\lambda})$--after allowing some additional ramification and enlarging the subgroup $H_1 \cap C^{\vee}$ of $C^{\vee}$--we have to carry out this annihilation in a way that is independent of $\lambda$, and moreover for fixed $\lambda$ independent of $\rho_{\lambda}$. Simultaneous annihilation of the $\mc{O}(\rho_{\lambda})$ using only a uniform, finite enlargement of the allowable ramification set and of the subgroup of $C^{\vee}$ in fact does not seem to be possible: we will as a first step have to define modified versions of these obstruction classes that take into account the Hodge numbers of $\psi$.
\end{rmk}
Before proceeding, we reinterpret the obstruction $\mc{O}(\rho_{\lambda})$ (we will only use the local version of what follows; in particular, the arguments of the present section depend only on the local version of Tate's theorem, which is an almost immediate consequence of local duality):
\begin{lemma}\label{Odescribed}
Let $v$ be a finite place of $F$, and suppose that $\tilde{\rho}_{\lambda} \colon \gal{F_v} \to \tH(\oE_{\lambda})$ is any continuous homomorphism lifting $\rho_{\lambda}|_{\gal{F_v}}$. Then $\mc{O}(\rho_{\lambda})|_{\gal{F_v}}$ is equal to the inverse of $\mc{O}(\nu \circ \tilde{\rho}_{\lambda})$, the obstruction associated to lifting $\nu \circ \tilde{\rho}_{\lambda} \colon \gal{F_v} \to \tZ^\vee(\oE_{\lambda})$ to $C^{\vee}$. (The same holds if we replace $\gal{F_v}$ by $\gal{F}$, but we do not require this.)
\end{lemma}
\proof
Before beginning the proof proper, we make precise our convention for co-boundary maps: the \textit{inverse} appearing in the conclusion of the lemma is crucial, and it is easy to get confused if one is not careful with the definitions. Let $\Gamma$ be a group and $M$ a (for simplicity) trivial $\Gamma$-module. For a function $\alpha \colon \Gamma^n \to M$, set 
\[
\delta(\alpha)(g_1, \ldots, g_{n+1})= \alpha(g_2, \ldots, g_n)+ \sum_{i=1}^n (-1)^i \alpha(g_1, \ldots, g_i g_{i+1}, \ldots, g_{n+1})+ (-1)^{n+1} \alpha(g_1, \ldots, g_n).
\]
For $n=1$, this says $\delta(\alpha)(g,h)= \alpha(h) \alpha(gh)^{-1} \alpha(g)$, and in the situation considered above $(g,h) \mapsto \rho'_{\lambda}(gh) \rho'_{\lambda}(h)^{-1} \rho'_{\lambda}(g)^{-1}$ is in fact a 2-cocycle.

Tate's theorem implies that, for large enough $M$, the image of $\mc{O}(\rho_{\lambda})$ in $H^2(\gal{F_v}, C^{\vee}[M])$ vanishes, i.e. $\mc{O}(\rho_{\lambda})= \delta (\phi)$ for some $\phi \colon \gal{F_v} \to C^{\vee}[M]$. The product $\rho'_{\lambda} \cdot \phi$ is then a homomorphism $\gal{F_v} \to (H_1 \cdot C^{\vee}[M])(\oE_{\lambda})$ lifting $\rho_{\lambda}$; we set $\tilde{\rho}_{\lambda}= \rho'_{\lambda} \cdot \phi$. Clearly $\nu \circ \tilde{\rho}_{\lambda}= \nu \circ \phi$, and then $\mc{O}(\nu \circ \tilde{\rho}_{\lambda})$ is (tautologically) represented by the co-cycle $(g, h) \mapsto \phi(gh) \phi(h)^{-1} \phi(g)^{-1}$, i.e. by $\delta(\phi)^{-1}= \mc{O}(\rho_{\lambda})^{-1} \in Z^2(\gal{F}, H_1 \cap C^{\vee})$.\footnote{Note that $\phi$ is valued in $C^{\vee}[M]$, not $H_1 \cap C^{\vee}$, so $\delta(\phi)$ need not be a co-boundary in $Z^2(\gal{F_v}, H_1 \cap C^{\vee})$.} This proves the claim for our particular lift $\tilde{\rho}_{\lambda}$, but any other lift $\tilde{\rho}^1_{\lambda}$ gives rise to the same obstruction $\mc{O}(\nu \circ \tilde{\rho}^1_{\lambda})$. (The global claim holds for the same reasons, if we admit the global version of Tate's theorem.)
\endproof
To address the difficulty indicated in Remark \ref{thepoint}, we begin by using the abelian representations coming from $\psi$ to construct a second obstruction class. Namely, consider the realization $\psi_{\iota_{\lambda}}$, which for notational simplicity from now on we simply denote by $\psi_{\lambda}$. The automorphic representation $\psi$ is unramified outside the finite set of places $T$ of $F$, so $\psi_{\lambda}$ is a geometric representation $\gal{F, T \cup S_{\lambda}} \to \tZ^\vee(\oE_{\lambda})$, which has good reduction outside $T$ (i.e. is crystalline at primes of $S_{\lambda}$ not in $T$). Via the isogeny $C^{\vee} \to \tZ^\vee$, we can then form a cohomology class measuring the obstruction to lifting $\psi_{\lambda}$ to $C^{\vee}$: let $\psi'_{\lambda}$ denote a topological lift $\gal{F, T \cup S_{\lambda}} \to C^{\vee}(\oE_{\lambda})$, defining as before a cohomology class
\[
\mc{O}(\psi_{\lambda}) \in H^2(\gal{F, T \cup S_{\lambda}}, H_1 \cap C^{\vee}).
\]
We can in turn define (via inflation) a cohomology class 
\[
\mc{O}(\rho_{\lambda}, \psi_{\lambda})=\mc{O}(\rho_{\lambda}) \cdot \mc{O}(\psi_{\lambda}) \in H^2(\gal{F, S \cup T \cup S_{\lambda}}, H_1 \cap C^{\vee}),
\]
which is represented by the 2-cocycle (recall that $C^{\vee}$ is central in $\tH$)
\[
(g, h) \mapsto (\rho'_{\lambda}\cdot \psi'_{\lambda})(gh) (\rho'_{\lambda}\cdot \psi'_{\lambda})(h)^{-1}(\rho'_{\lambda}\cdot \psi'_{\lambda})(g)^{-1}.
\]
(Note, however, that the function $g \mapsto (\rho'_{\lambda}\cdot \psi'_{\lambda})(g)$ is valued in $\tH$, not in $H_1$.)

We need one more lemma before getting to the crucial local result (Lemma \ref{loctriv} below)
\begin{lemma}\label{crislift}
For all places $v \in S_{\lambda}$, and for any choice of embedding $\iota_{\lambda} \colon \oE \into \oE_{\lambda}$, there exists a de Rham lift 
\[
\xymatrix{
& \tH(\oE_{\lambda}) \ar[d] \\
\gal{F_v} \ar[r]_-{\rho_{\lambda}} \ar@{-->}^{\tilde{\rho}_{\lambda}}[ur] & H(\oE_{\lambda})
}
\]
of $\rho_{\lambda}|_{\gal{F_v}}$ such that for all embeddings $\tau_{\lambda} \colon F_v \into \oE_{\lambda}$, the $\tau_{\lambda}$-labeled Hodge-Tate co-character of $\tilde{\rho}_{\lambda}$ is (conjugate to) $\wt{\mu_{\tau}} \otimes_{\oE, \iota_{\lambda}} \oE_{\lambda}$, where $\tau \colon F \into \oE$ is defined by the diagram
\[
\xymatrix{
F_v \ar[r]^{\tau_{\lambda}} & \oE_{\lambda} \\
F \ar[u] \ar[r]^{\tau} & \oE \ar[u]^{\iota_{\lambda}}.
}
\]
Moreover, if $\rho_{\lambda}$ is crystalline, then $\tilde{\rho}_{\lambda}$ may be taken to be crystalline.
\end{lemma}
\proof
For each $\tau_{\lambda} \colon F_v \into \oE_{\lambda}$, set for notational simplicity $\wt{\mu_{\tau_{\lambda}}}= \wt{\mu_{\tau}} \otimes_{\oE, \iota_{\lambda}} \oE_{\lambda}$, where $\tau$ is determined as in the diagram, and where $\wt{\mu_{\tau}}$ is the lift of $\mu_{\tau}$ we have fixed above. The proof of \cite[Corollary 3.2.12]{stp:variationsarxiv} shows that for any collection of co-characters lifting the Hodge co-characters of $\rho_{\lambda}$, and in particular for our $\wt{\mu_{\tau_{\lambda}}}$, there exists a Hodge-Tate lift $\tilde{\rho}_{\lambda} \colon \gal{F_v} \to \tH(\oE_{\lambda})$ whose $\tau_{\lambda}$-labeled Hodge-Tate co-character is $\wt{\mu_{\tau_{\lambda}}}$. Now consider the isogeny lifting problem
\[
\xymatrix{
&& \tH(\oE_{\lambda}) \ar[d] \\
\gal{F_v} \ar[rr]_-{\left(\rho_{\lambda}, \nu(\tilde{\rho}_{\lambda})\right)} \ar@{-->}[urr] && H(\oE_{\lambda}) \times \tZ^{\vee}(\oE_{\lambda}).
}
\]
Since $(\rho_{\lambda}, \nu(\tilde{\rho}_{\lambda}))$ admits a Hodge-Tate lift (namely, $\tilde{\rho}_{\lambda}$), and is itself de Rham ($\rho_{\lambda}$ is de Rham by assumption, and any abelian Hodge-Tate representation is de Rham), we can apply \cite[Corollary 6.7]{conrad:dualGW} to deduce the existence of a de Rham lift $\tilde{\rho}'_{\lambda}$, which clearly has the same Hodge-Tate co-characters as $\tilde{\rho}_{\lambda}$, since they differ by a finite-order twist. If we further assume $\rho_{\lambda}$ is crystalline, then we need only a minor modification to this argument: some power $\nu(\tilde{\rho}_{\lambda})^d$ is crystalline, so if we instead consider the problem of lifting the crystalline representation $(\rho_{\lambda}, [d]\nu(\tilde{\rho}_{\lambda}))$ through the composite isogeny 
\[
\tH \to H \times \tZ^\vee \xrightarrow{{\id \times [d]}} H \times \tZ^\vee,
\]
then again \cite[Corollary 6.7]{conrad:dualGW} applies to produce a crystalline lift of $\rho_{\lambda}$ with the desired Hodge-Tate co-characters.
\endproof

 Here is the key lemma:
\begin{lemma}\label{loctriv}
For any place $v \in S_{\lambda}$ not belonging to the finite set $S \cup T$, the restriction $\mc{O}(\rho_{\lambda}, \psi_{\lambda})|_{\gal{F_v}}$ is trivial.
\end{lemma}
\proof
Under the assumption on $v$, both $\rho_{\lambda}$ and $\psi_{\lambda}$ are crystalline at $v$; assume this from now on. Lemma \ref{crislift} above shows that $\rho_{\lambda}|_{\gal{F_v}}$ admits a crystalline lift $\tilde{\rho}_{\lambda} \colon \gal{F_v} \to \tH(\oE_{\lambda})$ such that $\nu \circ \tilde{\rho}_{\lambda}$ has the same (labeled) Hodge-Tate co-characters as $\psi_{\lambda}|_{\gal{F_v}}$. Since they are both crystalline, it follows (see \cite[3.9.7 Corollary]{chai-conrad-oort:cm}) that $(\nu \circ \tilde{\rho}_{\lambda}) \cdot \psi_{\lambda}^{-1}|_{\gal{F_v}}$ is unramified; this is an elaboration of the familiar fact that a crystalline character whose Hodge-Tate weights are zero must be unramified. In particular, replacing the initial lift $\tilde{\rho}_{\lambda}|_{\gal{F_v}}$ by an unramified twist, we may assume $\nu \circ \tilde{\rho}_{\lambda}= \psi_{\lambda}$ as homomorphisms $\gal{F_v} \to \tZ^\vee(\oE_{\lambda})$. But recall that Lemma \ref{Odescribed} implies that $\mc{O}(\rho_{\lambda})=\mc{O}(\nu \circ \tilde{\rho}_{\lambda})^{-1}$, so we deduce that $\mc{O}(\rho_{\lambda}) \cdot \mc{O}(\psi_{\lambda})|_{\gal{F_v}}$ is trivial.
\endproof

Since the set of places $S \cup T$ is finite, by the local version of Tate's theorem, the vanishing of $H^2(\gal{F_v}, \Q/\Z)$ for all places $v$ of $F$, we may enlarge $H_1 \cup C^{\vee}$ to some $C^{\vee}[m]$ inside the torus $C^{\vee}$ so as to kill the image of $H^2(\gal{F_v}, H_1 \cap C^{\vee}) \to H^2(\gal{F_v}, C^{\vee}[m])$ for all $v \in S \cup T$. (We emphasize that $m$ only depends on the set of places $S \cup T$ of $F$ and the finite group $H_1 \cap C^{\vee}$.) It follows then from Lemma \ref{loctriv} that if $\lambda$ does not belong to $S \cup T$, then $\mc{O}(\rho_{\lambda}, \psi_{\lambda})$ in fact belongs to 
\[
\Sh^2_{S \cup T \cup S_{\lambda}}(F, C^{\vee}[m])= \ker \left( H^2(\gal{F, S \cup T \cup S_{\lambda}}, C^{\vee}[m]) \to \bigoplus_{v \in S \cup T \cup S_{\lambda}} H^2(\gal{F_v}, C^{\vee}[m]) \right).
\]
We can moreover guarantee that this holds regardless of $\lambda$ by an additional finite enlargement of $m$ (since the number of exceptional $\lambda$ is finite). Moreover, by (if necessary) including the primes dividing $\#(H_1 \cap C^{\vee})$ in $S \cup T$, we can assume that $m$ is divisible only by primes in $S \cup T$. (Note that inflation to allow additional primes of ramification still has image in the corresponding Shafarevich-Tate group, since $\gal{F_v}/I_{F_v}$ has cohomological dimension one for all finite places $v$.) Thus, after these uniform enlargements of $m$ and $S \cup T$ (which we do not reflect in the notation), we have $\mc{O}(\rho_{\lambda}, \psi_{\lambda}) \in \Sh^2_{S \cup T \cup S_{\lambda}}(F, C^{\vee}[m])$.

We are now in a position to apply global duality to analyze the cohomology group $\Sh^2_{S \cup T \cup S_{\lambda}}(F, C^{\vee}[m])$. We will need, however, to allow still more primes of ramification in order to kill the class $\mc{O}(\rho_{\lambda}, \psi_{\lambda})$; the following crucial lemma allows us to do this in a way that does not depend on $\lambda$, but before stating the lemma, we have to recall the Gr\"{u}nwald-Wang theorem (somewhat specialized):
\begin{thm}[Gr\"{u}nwald-Wang; see Theorem X.1 of \cite{artin-tate:cft}]
Let $F$ be a number field, and let $m$ be a positive integer. Then an element $x \in F^\times$ belongs to $(F^\times)^m$ if and only if $x$ is in $(F_v^{\times})^m$ for all places $v$ of $F$, except when all three of the following conditions, to be referred to as the \textrm{special case}, hold for the pair $(F, m)$:
\begin{itemize}
\item Let $s_F$ denote the largest integer $r$ such that $\eta_r= \zeta_{2^r} + \zeta_{2^r}^{-1}$ is an element of $F$ (here $\zeta_{2^r}$ denotes a primitive $2^r$-th root of unity). Then $-1$, $2+ \eta_{s_F}$, and $-(2+\eta_{s_F})$ are non-squares in $F$.
\item $\mr{ord}_2(m)>s_F$.
\item The set of $2$-adic places of $F$ at which $-1$, $2+ \eta_{s_F}$, and $-(2+\eta_{s_F})$ are non-squares in $F$ is empty.
\end{itemize}
In the special case, the element $(2+ \eta_{s_F})^{m/2}$ is the unique (up to $(F^{\times})^m$-multiple) counter-example to the local-global principle for $m^{th}$-powers in $F^\times$.
\end{thm}
Here is the lemma:
\begin{lemma}
Recall that $S \cup T$ is a fixed finite set of places of the number field $F$, and that $m$ is a fixed integer. Let $V$ be a finite set of finite places of $F$ such that
\begin{itemize}
\item all elements of $V$ are unramified in $F(\mu_m)$;
\item the places of $F(\mu_m)$ lying above $V$ generate the class group of $F(\mu_m)$; and
\item every element of $\Gal(F(\mu_m)/F)$ is equal to a (geometric, say) frobenius element at $v$ for some $v \in V$.
\end{itemize}
Then for all places $\lambda$ of $E$ we can deduce: 
\begin{enumerate}
\item If $(F, m)$ is not in the Gr\"{u}nwald-Wang special case, $\Sh^2_{S \cup T \cup V \cup S_{\lambda}}(F, C^{\vee}[m])$ is trivial.
\item If $(F, m)$ is in the Gr\"{u}nwald-Wang special case, then the image of the canonical map
\[
\Sh^2_{S \cup T \cup V \cup S_{\lambda}}(F, C^{\vee}[m]) \to \Sh^2_{S \cup T \cup V \cup S_{\lambda}}(F, C^{\vee}[2m])
\]
is trivial.
\end{enumerate}
\end{lemma}
\proof
First note that such sets $V$ exist, by finiteness of the class number and the \v{C}ebotarev density theorem. Since (all places of $F$ above) the primes dividing $m$ are contained in $S \cup T$, an application of Poitou-Tate duality immediately reduces us to showing (as a Galois module, $C^{\vee}[m]$ is $\dim (C^{\vee})$ copies of $\Z/m$)
\begin{enumerate}
\item if $(F, m)$ is not in the Gr\"{u}nwald-Wang special case, then $\Sh^1_{S \cup T \cup V \cup S_{\lambda}}(F, \mu_m)=0$.
\item if $(F, m)$ is in the special case, then the map 
\[
\Sh^1_{S \cup T \cup V \cup S_{\lambda}}(F, \mu_{2m}) \to \Sh^1_{S \cup T \cup V \cup S_{\lambda}}(F, \mu_m)
\]
induced by $\mu_{2m} \xrightarrow{2} \mu_m$ is zero.
\end{enumerate}  
We first restrict to $\gal{F(\mu_m), S \cup T \cup V \cup S_{\lambda}}$ (note that this is actually restriction to a subgroup, since $F(\mu_m)/F$ is ramified only at primes in $S \cup T$), obtaining an element of $\Sh^1_{S \cup T \cup V \cup S_{\lambda}}(F(\mu_m), \mu_m)$. After this restriction, as we will see, the Gr\"{u}nwald-Wang theorem does not intervene.

To lighten the notation in the rest of the proof, we set $L= F(\mu_m)$ and $Q_{\lambda}= S \cup T \cup V \cup S_{\lambda}$. We also refer the reader to the notation established in \S \ref{notation}. Recall that $F(Q_{\lambda})$ denotes the maximal extension of $F$ inside $\overline{F}$ that is unramified outside $Q_{\lambda}$; it contains $L$. Let $\mc{O}_{F(Q_{\lambda})}$ denote the ring of $Q_{\lambda}$-integers in $F(Q_{\lambda})$ (i.e., the elements of $F(Q_{\lambda})$ that are integral outside of places above $Q_{\lambda}$). We then have an exact (Kummer theory) sequence
\[
1 \to \mu_m \to \mc{O}_{F(Q_{\lambda})}^{\times} \xrightarrow{m} \mc{O}_{F(Q_{\lambda})}^{\times} \to 1,
\]
and the corresponding long exact sequence in $\gal{L, Q_{\lambda}}$-cohomology yields an isomorphism
\[
\mc{O}_{L}[\frac{1}{Q_{\lambda}}]^{\times}/ \left( \mc{O}_{L}[\frac{1}{Q_{\lambda}}]^{\times} \right)^m \xrightarrow{\sim} H^1(\gal{L, Q_{\lambda}}, \mu_m);
\]
critically, surjectivity here follows from the vanishing of $H^1(\gal{L, Q_{\lambda}}, \mc{O}_{F(Q_{\lambda})}^{\times})$, which itself is a consequence of the natural isomorphism $Cl_{Q_{\lambda}}(L) \cong H^1(\gal{L, Q_{\lambda}}, \mc{O}_{F(Q_{\lambda})}^{\times})$ (see \cite[Proposition 8.3.11(ii)]{neukirch:cohnum}) and our assumption that $V$ (hence $Q_{\lambda}$) generates the class group of $L$. Restricting the Kummer theory isomorphism to classes that are locally trivial at each place of $Q_{\lambda}$, we also obtain the isomorphism
\[
\left( \mc{O}_L[\frac{1}{Q_{\lambda}}]^{\times} \cap (L_{Q_{\lambda}}^{\times})^m \right)/ \left( \mc{O}_{L}[\frac{1}{Q_{\lambda}}]^{\times} \right)^m \xrightarrow{\sim} \Sh^1_{Q_{\lambda}}(L, \mu_m).
\]
We claim these groups are trivial. Indeed, let $\alpha$ be an element of $\mc{O}_L[\frac{1}{Q_{\lambda}}]^{\times} \cap (L_{Q_{\lambda}}^{\times})^m$, and consider the (abelian) extension $L(\alpha^{1/m})/L$. Global class field theory yields the reciprocity isomorphism 
\[
\mathbb{A}_L^{\times}/ \left( L^{\times} N_{L(\alpha^{1/m})/L}(\mathbb{A}_{L(\alpha^{1/m})}^{\times}) \right) \xrightarrow{\sim} \Gal(L(\alpha^{1/m})/L),
\]
but by assumption the source of this map admits a surjection
\[
\mathbb{A}_L^{\times}/ \left( L^{\times} L_{\infty}^{\times} L_{Q_{\lambda}}^{\times} \prod_{w \not \in Q_{\lambda}} \mc{O}_{L_w}^{\times} \right) \onto \mathbb{A}_L^{\times}/\left( L^{\times} N_{L(\alpha^{1/m})/L}(\mathbb{A}_{L(\alpha^{1/m})}^{\times}) \right).
\]
(At unramified places, the image of the norm map contains the local units; and at places in $Q_{\lambda}$, $L(\alpha^{1/m})/L$ is split.) By assumption ($Cl_{Q_{\lambda}}(L)=0$), the source of this surjection is trivial, so $L(\alpha^{1/m})=L$, and we deduce that $\Sh^1_{Q_{\lambda}}(L, \mu_m)=0$.

It follows that inflation identifies $\Sh^1_{Q_{\lambda}}(F, \mu_m)$ with the classes in $H^1(\Gal(L/F), \mu_m)$ that are trivial upon restriction to $Q_{\lambda}$. Since every element of $\Gal(L/F)$ is a frobenius element at some prime in $V \subset Q_{\lambda}$, $\Sh^1_{Q_{\lambda}}(F, \mu_m)$ is actually equal to the set of \textit{everywhere} locally trivial classes 
\[
\Sh^1_{|F|}(F, \mu_m):= \ker\left( H^1(\gal{F}, \mu_m) \to \prod_{v \in |F|} H^1(\gal{F_v}, \mu_m) \right),
\]
where $|F|$ denotes the set of all places of $F$. This is precisely the subject of the Gr\"{u}nwald-Wang theorem, and it is zero if $(F, m)$ is not in the special case. Thus, we need only consider the possibility that $(F, m)$ is in the special case, where $\Sh^1_{|F|}(F, \mu_m)$ has order two, and a representative of the non-trivial class is the (image under the Kummer map of the) element $(2+ \eta_{s_F})^{m/2}$ of $(F^\times)^{m/2}$. This description holds regardless of $m$, so in particular the non-trivial class of $\Sh^1_{|F|}(F, \mu_{2m})$ is represented by $(2+ \eta_{s_F})^{m}$. Its image under $\mu_{2m} \xrightarrow{2} \mu_m$, which via Kummer theory is induced by the identity map $F^\times \to F^\times$, is again $(2+ \eta_{s_F})^m$, which is now visibly an $m^{th}$ power, completing the proof.\footnote{In concrete terms, this says that if an element of $F^\times$ is everywhere locally a $(2m)^{th}$-power, then it is globally an $m^{th}$-power.}
\endproof
We summarize our conclusion, noting that the value of $m$ in the following corollary may be $2m$ in the earlier notation:
\begin{cor}\label{killedcoh}
There is an integer $m$ and a finite set of places $Q \supset S \cup T$, both independent of $\lambda$ and of the choice of $\rho_{\lambda}$ (having good reduction outside $S$ and prescribed Hodge-Tate co-characters), such that the image of $\mc{O}(\rho_{\lambda}, \psi_{\lambda})$ in $H^2(\gal{F, Q \cup S_{\lambda}}, C^{\vee}[m])$ is zero.
\end{cor}
Before proceeding, it is worth noting that the argument just given yields a novel proof of the global version of Tate's vanishing theorem (taking as input the much easier local theorem); it is also a stronger proof, yielding an explicit upper-bound on how much ramification has to be allowed, and how much the coefficients need to be enlarged, in order to kill a cohomology class in $H^2(\gal{F, V}, \Z/N)$ for some finite set of places $V$ and integer $N$:
\begin{cor}\label{tate}
Let $V$ be a finite set of places of $F$, and let $N$ be an integer. Then the image of $H^2(\gal{F, V}, \Z/N)$ in $H^2(\gal{F, V \cup W}, \Z/2NM)$ is trivial where:
\begin{itemize}
\item $M$ is large enough that for all $v \in V$, the image of $H^2(\gal{F_v}, \Z/N) \to H^2(\gal{F_v}, \Z/NM)$ is zero;\footnote{This is easy to make explicit, using local duality, in terms of $\mu_{\infty}(F_v)$.} and
\item once $M$ is fixed as above, $W$ is large enough that 
\begin{itemize}
\item $V \cup W$ contains (all places above) $2NM$, 
\item $Cl_{V \cup W}(F(\mu_{NM}))=0$, and
\item each element of $\Gal(F(\mu_{NM})/F)$ is equal to a frobenius element at $w$ for some $w \in W$.
\end{itemize}
\end{itemize}
(The factors of two are only necessary in the Gr\"{u}nwald-Wang special case) In particular, $H^2(\gal{F}, \Q/\Z)=0$.
\end{cor}
\begin{rmk}
A different proof of Tate's theorem (without arithmetic duality theorems, but instead relying on a finer study of Hecke characters of $F$) is given in \cite[\S 6.5]{serre:DSsurvey}. There Serre remarks that Tate originally proved the vanishing theorem using global duality, but further assuming Leopoldt's conjecture; we have of course circumvented Leopoldt here. 
\end{rmk}
Now we return to the conclusion of Corollary \ref{killedcoh}. Let $b_{\lambda} \colon \gal{F, Q \cup S_{\lambda}} \to C^{\vee}[m]$ be a cochain trivializing $\mc{O}(\rho_{\lambda}, \psi_{\lambda})$. Then
\[
\tilde{\rho}_{\lambda}= \rho'_{\lambda} \cdot \psi'_{\lambda} \cdot b_{\lambda} \colon \gal{F, Q \cup S_{\lambda}} \to \tH(\oE_{\lambda})
\]
is a homomorphism lifting $\rho_{\lambda}$. We claim that $\tilde{\rho}_{\lambda}$ is moreover de Rham at all places in $S_{\lambda}$. To see this, note that under the isogeny $\tH \to H \times \tZ^\vee$, $\tilde{\rho}_{\lambda}$ pushes forward to $(\rho_{\lambda}, \psi_{\lambda} \nu(b_{\lambda}))$, the second coordinate being a finite-order twist of $\psi_{\lambda}$ (and in particular, de Rham). But now we can invoke the local results of Wintenberger (\cite[\S 1]{wintenberger:relevement}) and Conrad (\cite[Theorem 6.2]{conrad:dualGW}), asserting that a lift of a de Rham representation through an isogeny is de Rham if and only if the Hodge-Tate co-character lifts through the isogeny (which is obviously the case here, as $\psi$ was constructed to ensure this). 

Finally, we can refine this to the statement that $\rho_{\lambda}$ admits a geometric lift that is moreover crystalline at all places of $S_{\lambda}$, provided $S_{\lambda}$ does not intersect a certain finite set of primes that is independent of $\lambda$ and $\rho_{\lambda}$ (but somewhat larger than the set $Q$ we have thus far constructed):
\begin{thm}[Theorem \ref{mainintro}]\label{main}
Let $F$ be a totally imaginary number field, and let $S$ be a finite set of places of $F$ containing the archimedean places. Fix a set of co-characters $\{\mu_{\tau}\}_{\tau \colon F \into \oE}$ as in part (2) of Definition \ref{ramcomp}. Then there exists a finite set of places $P \supset S$ such that any geometric representation $\rho_{\lambda} \colon \gal{F, S \cup S_{\lambda}} \to H(\oE_{\lambda})$ having good reduction outside $S$, and whose Hodge-Tate co-characters arise from the set $\{\mu_{\tau}\}_{\tau \colon F \into \oE}$ via an embedding $\iota_{\lambda} \colon \oE \into \oE_{\lambda}$, admits a geometric lift $\tilde{\rho}_{\lambda} \colon \gal{F, P \cup S_{\lambda}} \to \tH(\oE_{\lambda})$ having good reduction outside $P$.

In particular, if $\{\rho_{\lambda} \colon \gal{F, S \cup S_{\lambda}} \to H(\oE_{\lambda})\}_{\lambda}$ is a ramification-compatible system, then there exist a finite set of places $P \supset S$ and lifts $\tilde{\rho}_{\lambda} \colon \gal{F, P \cup S_{\lambda}} \to \tH(\oE_{\lambda})$ such that $\{\tilde{\rho}_{\lambda}\}_{\lambda}$ is a ramification-compatible system.
\end{thm}
\proof
We resume the above discussion. So far we have a constructed geometric lifts
\[
\tilde{\rho}_{\lambda} \colon \gal{F, Q \cup S_{\lambda}} \to \tH(\oE_{\lambda}),
\]
where $Q$ contains $S \cup T$ and whatever other additional primes are needed for the conclusion of Corollary \ref{killedcoh}. The only remaining task is to show that for some (independent of $\rho_{\lambda}$) set $P$, we can modify the initial lift (by a finite-order twist) to guarantee that it has good reduction outside $P$. Under the isogeny $\tH \to H \times \tZ^\vee/\nu(C^{\vee}[m])$, $\tilde{\rho}_{\lambda}$ pushes forward to 
\[
\tau_{\lambda} := \left(\rho_{\lambda}, \psi_{\lambda}\mod {\nu (C^{\vee}[m])}\right),
\] 
which is crystalline for all $v$ in $S_{\lambda}$ but not in $S \cup T$. For all $v \in S_{\lambda} \setminus \left(S_{\lambda} \cap (S \cup T) \right)$, $\tilde{\rho}_{\lambda}|_{\gal{F_v}}$ is of course a de Rham lift of $\tau_{\lambda}$, so \cite[Theorem 6.2, Corollary 6.7]{conrad:dualGW} (building on \cite{wintenberger:relevement}) shows that $\tau_{\lambda}|_{\gal{F_v}}$ admits \textit{some} crystalline lift $\tilde{\tau}_{\lambda, v} \colon \gal{F_v} \to \tH(\oE_{\lambda})$, and therefore there are finite-order characters $\chi_{\lambda, v} \colon \gal{F_v} \to C^{\vee}[m]$ such that each $\tilde{\rho}_{\lambda}|_{\gal{F_v}} \cdot \chi_{\lambda, v}$ is crystalline. We wish to glue the inertial restrictions $\chi_{\lambda, v}|_{I_{F_v}}$ together into a global character, with an independent-of-$\lambda$ control on the ramification. The cokernel of the restriction map
\begin{equation}\label{restrict}
\Hom(\gal{F, Q \cup S_{\lambda}}, C^{\vee}[m]) \to \bigoplus_{v \in S_{\lambda}} \Hom(I_{F_v}, C^{\vee}[m])^{\gal{F_v}/I_{F_v}}
\end{equation}
may be non-trivial;\footnote{Of course we need only consider the cokernel of the map to the direct sum over $v \in S_{\lambda} \setminus (S_{\lambda} \cap (S \cup T))$; to lighten the notation we will work with all $v \in S_{\lambda}$, taking some arbitrary (eg, trivial) choice of $\chi_{\lambda, v}$ at any places in $(S \cup T) \cap S_{\lambda}$.} but we will show that any element of the cokernel is annihilated by appropriate enlargements of $Q$ and $m$.

By a theorem of Chevalley (the congruence subgroup property for $\mr{GL}_1$), there is an ideal $\mf{n}$ of $\mc{O}_F$ such that
\[
\left\{x \in \mc{O}_F^\times: x \equiv 1 \pmod {\mf{n}}\right\} \subseteq (\mc{O}_F^\times)^m.
\]
Let $R$ be the set of primes supporting $\mf{n}$ (note that $\mf{n}$ and $R$ are independent of $\rho_\lambda$!), and set 
\[
U_R= \left\{(x_v)_{v \in R} \in \prod_{v \in R} \mc{O}_{F_v}^\times: x_v \equiv 1 \pmod{\mf{n}} \quad \text{for all $v \in R$}\right\}.
\]
Then whenever $S_{\lambda} \cap R= \emptyset$ (so, excluding a finite number of bad $\lambda$), consider the character (here and in what follows, we suppress the class field theory identifications)
\[
(\chi_{\lambda, v})_{v \in S_{\lambda}} \times 1 \times 1 \times 1 \colon\prod_{v \in S_{\lambda}} \mc{O}_{F_v}^\times \times \prod_{v \in R} U_R \times \prod_{v \not \in R \cup S_{\lambda}} \mc{O}_{F_v}^\times \times F_{\infty}^\times \to C^{\vee}[m],
\]
which extends by 1 to a character 
\[
\left(\prod_{v \in S_{\lambda}} \mc{O}_{F_v}^\times \times \prod_{v \in R} U_R \times \prod_{v \not \in R \cup S_{\lambda}} \mc{O}_{F_v}^\times \times F_{\infty}^\times \right) \cdot F^\times \to C^{\vee}[m]
\]
(an element of the intersection is a global unit congruent to 1 modulo $\mf{n}$, hence is contained in $(\mc{O}_F^\times)^m$, where $\chi_{\lambda, v}$ is obviously trivial). We can then extend from this finite-index subgroup of $\mathbb{A}_F^\times$ to a character $\chi_{\lambda} \colon \mathbb{A}_F^\times/F^{\times} \to \mu_{\infty}(C^{\vee})$. In fact, we see that $\chi_{\lambda}$ can be chosen to be valued in $C^{\vee}[M]$ for $M$ sufficiently large but independent of $\lambda$: $M$ can be quantified in terms of the generalized class group of level $U_R$, but the details do not concern us.

Replacing $\tilde{\rho}_{\lambda}$ by its finite-order twist 
\[
\tilde{\rho}_{\lambda} \cdot \chi_{\lambda} \colon \gal{F, Q \cup R \cup S_{\lambda}} \to \tH(\oE_{\lambda}),
\]
we have achieved geometric lifts of $\rho_{\lambda}$ with compatible Hodge-Tate co-characters, and which are crystalline at all places in $S_{\lambda}$ outside of $R \cup S \cup T$.
\endproof
\begin{rmk}
Contrast the final step \cite[\S 2.3.5]{wintenberger:relevement} of Wintenberger's main theorem, where to ensure crystallinity of the lifts he makes a further finite base-change on $F$ (having already made several such in order to show lifts exist, as is necessary in his isogeny set-up), adding appropriate roots of unity and then passing to a Hilbert class field to kill a cokernel analogous to that of Equation (\ref{restrict}). As elsewhere, our argument is orthogonal to Wintenberger's, in allowing additional ramification and larger coefficients, rather than passing to a finite extension of $F$.
\end{rmk}
We now deduce some corollaries on finding lifts of ramification-compatible systems whose ``similitude characters" (determinant, Clifford norm, etc.) form strongly-compatible systems, in the sense that at all finite places their associated Weil group representations are isomorphic (see, eg, \cite[\S 5.1]{blggt:potaut}, where these are called \textit{strictly} compatible). As with Theorem \ref{main} and the preceding results, we show a somewhat stronger finiteness result, which applies to all representations with good reduction outside a fixed finite set $S$. These corollaries will follow from the above results and the Hermite-Minkowski finiteness theorem. 
\begin{cor}\label{charcomp}
Let $F$, $S$, and $\{\mu_\tau\}$ be as in the statement of Theorem \ref{main}, except now $F$ may be any number field. Then there exist a finite set of places $P \supset S$ and a finite extension $F'/F$ such that any geometric $\rho_{\lambda}$ with good reduction outside $S$, and with Hodge-Tate co-characters arising from $\{\mu_{\tau}\}$ via an embedding $\iota_{\lambda} \colon \oE \into \oE_{\lambda}$, admits a geometric lift $\tilde{\rho}_{\lambda} \colon \gal{F, P \cup S_{\lambda}} \to \tH(\oE_{\lambda})$ such that the restrictions
\[
\nu(\tilde{\rho}_{\lambda}) \colon \gal{F', P \cup S_{\lambda}} \to \tZ^\vee(\oE_{\lambda})
\]
are equal to the $\iota_\lambda$-adic realizations of the single (independent of $\lambda$ and $\rho_{\lambda}$) algebraic Hecke character $\psi$ of $\tZ(\mathbb{A}_{F})$.

In particular, let $\{ \rho_{\lambda} \colon \gal{F, S \cup S_{\lambda}} \to H(\oE_{\lambda}) \}_{\lambda}$ be a ramification-compatible system. Then there exist a ramification-compatible system of lifts $\{ \tilde{\rho}_{\lambda} \colon \gal{F, P \cup S_{\lambda}} \to \tH(\oE_{\lambda})\}_{\lambda}$, and a finite, independent-of-$\lambda$, extension $F'/F$ such that the restrictions
\[
\nu(\tilde{\rho}_{\lambda}) \colon \gal{F', P \cup S_{\lambda}} \to \tZ^\vee(\oE_{\lambda})
\]
form a strongly-compatible system.
\end{cor}
\proof
We may assume $F$ is totally imaginary. Consider the lifts $\tilde{\rho}_{\lambda} \colon \gal{F, P \cup S_{\lambda}} \to \tH(\oE_{\lambda})$ produced by Theorem \ref{main}. We write $\nu(\tilde{\rho}_{\lambda})= \psi_{\lambda} \cdot \eta_{\lambda}$, where $\eta_{\lambda} \colon \gal{F, P \cup S_{\lambda}} \to \tZ^\vee[M]$ is a finite-order character: the independent-of-$\lambda$ bound on the order was established within the proof of Theorem \ref{main}. Moreover, for all $v \in S_{\lambda} \setminus (S_{\lambda} \cap P)$, $\tilde{\rho}_{\lambda}$ and $\psi_{\lambda}$ are crystalline at $v$, so as long as $S_{\lambda} \cap P$ is empty, $\eta_{\lambda}$ factors through $\gal{F, P} \to \tZ^\vee[M]$ (we again use that a finite-order crystalline character is unramified). By the Hermite-Minkowski theorem, there are a finite number of such characters $\eta_{\lambda}$. For the finite number of bad $\lambda$ (at which $S_{\lambda} \cap P \neq \emptyset$), the same finiteness assertion holds. Thus, after a finite base-change $F'/F$, trivializing this finite collection of possible characters $\eta_{\lambda}$, we see that $\nu(\tilde{\rho}_{\lambda})|_{\gal{F', P \cup S_{\lambda}}}= \psi_{\lambda}|_{\gal{F', P \cup S_{\lambda}}}$ for all $\lambda$. The second part of the corollary follows since the $\lambda$-adic realizations of an abelian L-algebraic representation form a strongly compatible system (as is evident from the construction of $\psi_{\lambda}$, as in, eg \cite{serre:ladic}).
\endproof
We would like to upgrade this to a compatibility statement not just for the push-forwards $\nu(\tilde{\rho}_{\lambda})$, but for the full abelianizations $\tilde{\rho}^{\mr{ab}} \colon \gal{F, P \cup S_{\lambda}} \to \tH^{\mr{ab}}(\oE_{\lambda})$. Of course, such a result requires first (taking $\tH=H$) having the corresponding assertion for the abelianizations $\rho_{\lambda}^{\mr{ab}} \colon \gal{F, S \cup S_{\lambda}} \to H^{\mr{ab}}(\oE_{\lambda})$. Here, however, it is of course false without imposing further conditions on the system $\{\rho_{\lambda} \}_{\lambda}$ (see Remark \ref{coefficientrmk}). There are various conditions we might impose on the $\rho_{\lambda}$ to ensure (potential) compatibility of the $\rho_{\lambda}^{\mr{ab}}$. Perhaps most interesting is to restrict the coefficients of $\rho_{\lambda}^{\mr{ab}}$. To that end, we first prove a finiteness result for Galois characters:
\begin{lemma}\label{coefffin}
Let $F$ be a number field, and let $S$ be a finite set of places of $F$. Fix a finite extension $E'/E$ (inside $\oE$), a set $\{m_{\tau}\}_{\tau \colon F \into \oE}$ of integers satisfying the Hodge-symmetry condition of Definition \ref{ramcomp}, and an embedding $\iota_{\infty} \colon \oE \into \CC$. Then there exist a finite extension $F'/F$, and an algebraic Hecke character $\alpha$ of $\mathbb{A}_{F'}$, such that any geometric character $\omega_{\lambda} \colon \gal{F, S \cup S_{\lambda}} \to \oE_{\lambda}^\times$
\begin{itemize}
\item having good reduction outside $S$;
\item having labeled Hodge-Tate weights corresponding to $\{m_{\tau}\}$ via some embedding $\iota_{\lambda} \colon \oE \into \oE_{\lambda}$; and
\item for which $\omega_{\lambda}(fr_v)$ belongs to $(E')^\times$ for a density-one set of places $v$ of $F$; 
\end{itemize}
will upon restriction $\omega_{\lambda}|_{\gal{F'}}$ become isomorphic to the $\iota_{\lambda}$-adic realization of $\alpha$.
\end{lemma}
\proof
We may assume $F$ is totally imaginary. Invoking the Hodge symmetry hypothesis, we apply Lemma \ref{hecke} to produce an algebraic Hecke character $\alpha$ of $F$ whose archimedean components are given in terms of the $m_{\tau}$, exactly as in Lemma \ref{hecke} (with $\nu \circ \wt{\mu}_{\tau}=m_{\tau}$). Let $T$ denote the finite set of ramified places of $\alpha$, and let $\Q(\alpha)$ denote the field of coefficients of $\alpha$ (by definition the fixed field of all automorphisms of $\CC$ that preserve the non-archimedean component of $\alpha$; we will regard $\Q(\alpha)$ as a subfield of $\oE$ via our fixed $\iota_{\infty}$). Thus the $\iota_{\lambda}$-adic realizations $\alpha_{\lambda} \colon \gal{F, T \cup S_{\lambda}} \to \oE_{\lambda}^\times$ have labeled Hodge-Tate weights matching those of $\omega_{\lambda}$. Since $\Q(\alpha)$ contains the values $\alpha_{\lambda}(fr_v)$ for all $v \not \in T \cup S_{\lambda}$, and $\omega_{\lambda} \alpha_{\lambda}^{-1} \colon \gal{F, S \cup T \cup S_{\lambda}} \to \oE_{\lambda}^\times$ is finite-order (all of its Hodge-Tate weights are zero), we see that $(\omega_{\lambda} \alpha_{\lambda}^{-1})(fr_v)$ belongs to the finite (independent of $\lambda$) set $\mu_{\infty}(E' \Q(\alpha))$ for a density-one set of $v$. By \v{C}ebotarev, the character $\omega_{\lambda} \alpha_{\lambda}^{-1}$ takes all of its values in $\mu_{\infty}(E' \Q(\alpha))$. As long as $S_{\lambda} \cap (S \cup T)$ is empty, $\omega_{\lambda} \alpha^{-1}_{\lambda}$ is moreover unramified at $S_{\lambda}$ (because it is crystalline of finite-order), so as in Corollary \ref{charcomp}, there are (again by Hermite-Minkowski) a finite number of such characters $\omega_{\lambda}\alpha_{\lambda}^{-1}$. We deduce the existence of a single number field $F'$ over which $\omega_{\lambda}|_{\gal{F'}}= \alpha_{\lambda}|_{\gal{F'}}$, for any $\lambda$ and any $\omega_{\lambda}$ as in the statement of the lemma.
\endproof
We then deduce a potential compatibility statement for the full abelianizations $\tilde{\rho}_{\lambda}^{\mr{ab}}$:
\begin{cor}\label{abcomp}
 For simplicity, assume that $H^{\mr{ab}}$ is of multiplicative type. Let $F$, $S$, and $\{\mu_{\tau}\}$ be as in the statement of Theorem \ref{main}, except with $F$ now allowed to be any number field. Also fix a finite extension $E'$ of $E$. Then there exist a finite set of places $P \supset S$, a finite extension $F'/F$, and an algebraic Hecke character $\beta$ of the split group $\wt{D}$ over $F'$ whose dual group over $\oE$ is isomorphic to (and we fix such an isomorphism) $(H^{\mr{ab}})^0 \otimes_E \oE$, satisfying the following: if a geometric $\rho_{\lambda} \colon \gal{F, S \cup S_{\lambda}} \to H(\oE_{\lambda})$ 
\begin{itemize}
 \item has good reduction outside $S$;
 \item has Hodge-Tate co-characters arising from $\{\mu_{\tau}\}$ via $\iota_{\lambda} \colon \oE \into \oE_{\lambda}$;
 \item and admits, for some faithful representation $r$ of $H^{\mr{ab}}$, a density-one set of places $v$ of $F$ such that the characteristic polynomial $\mr{ch}(r \circ \tilde{\rho}_{\lambda}^{\mr{ab}})(fr_v)$ has coefficients in $E'$;
\end{itemize}
then there is a geometric lift $\tilde{\rho}_{\lambda} \colon \gal{F, P \cup S_{\lambda}} \to \tH(\oE_{\lambda})$ having good reduction outside $P$, such that the restriction $\tilde{\rho}_{\lambda}^{\mr{ab}} \colon \gal{F', P \cup S_{\lambda}} \to \tH^{\mr{ab}}(\oE_{\lambda})$ is equal to the $\iota_{\lambda}$-adic realization $\beta_{\lambda}$ of $\beta$.

In particular, let $\{ \rho_{\lambda} \colon \gal{F, S \cup S_{\lambda}} \to H(\oE_{\lambda})\}_{\lambda}$ be a ramification-compatible system, and assume that for some faithful representation $r$ of $H^{\mr{ab}}$, some number field $E'$, and for almost all $\lambda$, there is a density-one set of places $v$ of $F$ such that the characteristic polynomial $\mr{ch}(r \circ \tilde{\rho}_{\lambda}^{\mr{ab}})(fr_v)$ has coefficients in $E'$. Then there is a ramification-compatible system $\tilde{\rho}_{\lambda} \colon \gal{F, P \cup S_{\lambda}} \to \tH(\oE_{\lambda})$ lifting $\rho_{\lambda}$, and a finite extension $F'/F$ such that 
\[
\tilde{\rho}_{\lambda}^{\mr{ab}}|_{\gal{F'}} \colon \gal{F', P \cup S_{\lambda}} \to \tH^{\mr{ab}}(\oE_{\lambda})
\]
forms a strongly-compatible system.
\end{cor}
\proof
The proof follows familiar lines. 
Since $C^{\vee}$ is central, the abelianization $\tH^{\mr{ab}}$ is simply $\tH/H_1^{\mr{der}}$, so there is a natural map
\[
f \colon \tH^{\mr{ab}} \to \tH/H_1 \times H/\im(H_1^{\mr{der}})= \tZ^\vee \times H^{\mr{ab}},
\]
under which $\tilde{\rho}_{\lambda}^{\mr{ab}}$ pushes forward to $(\nu(\tilde{\rho}_{\lambda}), \rho_{\lambda}^{\mr{ab}})$. (We have chosen $\{\tilde{\rho}_{\lambda}\}_{\lambda}$ as in Theorem \ref{main} and Corollary \ref{charcomp}, of course.) First we claim that a conclusion analogous to that of the corollary holds for the pair $(\nu(\tilde{\rho}_{\lambda}), \rho_{\lambda}^{\mr{ab}})$, 
and certainly it suffices to check this independently for the two components. The assertion for $\nu(\tilde{\rho}_{\lambda})$ is Corollary \ref{charcomp}, and for $\rho_{\lambda}^{\mr{ab}}$ it follows easily from Lemma \ref{coefffin} (first reduce, by a finite-base-change, to the case where $H^{\mr{ab}}$ is connected, using the fact that $\pi_0(H^{\mr{ab}})$ is of course finite and independent of $\lambda$). Thus, letting $D$ denote a split torus whose dual group is identified with $(H^{\mr{ab}})^0$, there exists a finite extension $F_1/F$ such that $f(\tilde{\rho}_{\lambda}^{\mr{ab}})|_{\gal{F_1, P \cup S_{\lambda}}}$ is the $\iota_{\lambda}$-adic realization of a Hecke character (not depending on $\lambda$ or $\rho_{\lambda}$) of $\tZ \times D$. 

Now let $\beta$ be a Hecke character of $\wt{D}$ whose $\iota_{\lambda}$-adic realization $\beta_{\lambda} \colon \gal{F, T \cup S_{\lambda}} \to (\tH^{\mr{ab}})^0(\oE_{\lambda})$ has labeled Hodge-Tate co-characters matching those of $\tilde{\rho}_{\lambda}^{\mr{ab}}$. Since 
\[
f(\tilde{\rho}_{\lambda}^{\mr{ab}} \cdot \beta_{\lambda}^{-1})|_{\gal{F_1, P \cup T \cup S_{\lambda}}} 
\]
is automorphic (independently of $\lambda$, $\rho_{\lambda}$) of finite-order, it is trivial after a finite base-change $F_2/F_1$. Now observe that the kernel of $f$ is finite, so $(\tilde{\rho}_{\lambda}^{\mr{ab}} \cdot \beta_{\lambda}^{-1})|_{\gal{F_2, P \cup T \cup S_{\lambda}}}$ has finite-order, bounded only in terms of $\# \ker(f)$, and is crystalline away from $P \cup T$; as before, we find a further finite extension $F_3/F_2$ such that $(\tilde{\rho}_{\lambda}^{\mr{ab}} \cdot \beta_{\lambda}^{-1})|_{\gal{F_3}}=1$. The conclusion of the Corollary then holds with $F'=F_3$.
\endproof
\begin{rmk}\label{coefficientrmk}
\begin{itemize}
\item It does not suffice to ask for a fixed number field $E$ such that all $\rho_{\lambda}^{\mr{ab}}$ are valued in $H^{\mr{ab}}(E_{\lambda})$. For instance, taking $F= \Q$ and $S= \{p\}$, and for all $n$ choosing a prime $\ell_n \equiv 1 \pmod{\varphi(p^n)}$, we can define $\rho_{\ell_n} \colon \gal{\Q, \{p\}} \to \Q_{\ell_n}^\times$ as the composition of the mod $p^n$ cyclotomic character with an inclusion $(\Z/p^n\Z)^\times \into \mu_{\ell_n-1} \into \Q_{\ell_n}^\times$, and for all $\ell \not \in \{\ell_n\}_n$ we can take $\rho_{\ell}$ to be the trivial character. Then $\{\rho_{\ell} \colon \gal{\Q, \{p\}} \to \Q_\ell^\times\}_{\ell}$ is an abelian, ramification-compatible system that does not become a strongly-compatible system after any finite base-change.
\item Having only hypothesized ramification-compatibility for the $\{\rho_{\lambda}\}_{\lambda}$, we cannot hope for the stronger conclusion that the $\{\tilde{\rho}_{\lambda}^{\mr{ab}}\}_{\lambda}$ form a strongly-compatible system over $F$ itself.
\end{itemize}
\end{rmk}

\subsection{General multiplicative-type quotients}\label{multtype}
In fact, the argument of Theorem \ref{main} directly implies the main theorem of \cite{wintenberger:relevement}, as well as a generalization to lifting through quotients where the kernel is central of multiplicative type. We thus obtain an essentially different proof (and generalization) of Wintenberger's result. In this section, we briefly describe how this works.
\begin{cor}[Wintenberger]\label{wint}
Let $H_1 \to H$ be a central isogeny of linear algebraic groups over $E$, and let $S$ be a finite set of places of $F$. Then there exist a finite extension $F'/F$ and a finite set of places $P \supset S$ of $F$ such that any 
\[
\rho_{\lambda} \colon \gal{F, S \cup S_{\lambda}} \to H(E_{\lambda})
\]
having:
\begin{itemize}
\item good reduction outside $S$, and
\item labeled Hodge-Tate co-characters that lift to $H_1$,
\end{itemize}
lifts to a geometric representation $\rho'_{\lambda} \colon \gal{F', P \cup S_{\lambda}} \to H_1(E_{\lambda})$, which moreover has good reduction outside $P$.
\end{cor}
\proof
We begin by replacing $F$ by a finite extension $F_0$ such that image of $\rho_{\lambda}|_{\gal{F_0}}$ is contained in the image of $H_1(E_{\lambda}) \to H(E_{\lambda})$. That such an extension, depending only on $H_1 \to H$, $S$, and $F$, exists follows as in \cite[2.3.2]{wintenberger:relevement}, and we do not repeat the argument. We note, though, that making this construction in an independent-of-$\lambda$ fashion already uses liftability of the Hodge-Tate co-characters. (If we were not concerned with preserving $E_{\lambda}$-rationality of the lift, then we could skip this step.) It is then possible to build an obstruction class $\mc{O}(\rho_{\lambda}) \in H^2(\gal{F_0, S \cup S_{\lambda}}, \ker(H_1 \to H)(E_{\lambda}))$ via a topological lift $\rho_{\lambda}'$ to $H_1(E_{\lambda})$.

Embed $\ker(H_1 \to H) \otimes_E \oE$ into a torus $C^{\vee}$, and form the new group $\tH= (H_1 \times C^{\vee})/\ker(H_1 \to H)$, with the kernel embedded anti-diagonally. The surjection $\tH \to H$ now has kernel equal to a central torus $C^{\vee}$, and as before we let $\tZ^\vee$ be the (torus) quotient $\tH/H_1$. By hypothesis, we can lift the Hodge-Tate co-characters of $\rho_{\lambda}$ to $H_1$; when pushed-forward to $\tZ^\vee$, these lifts are of course trivial. Thus, in the notation of Lemma \ref{hecke}, we may take the trivial Hecke character $\psi=1$ of $\tZ(\mathbb{A}_{F_0})$. For topological lifts $\psi'_{\lambda}$ to $C^{\vee}(E_{\lambda})$ (as in Lemma \ref{Odescribed}) of the (trivial) $\lambda$-adic realizations $\psi_{\lambda}$, we may of course also take $\psi'_{\lambda}=1$. Theorem \ref{main} then produces a finite set of primes $P \supset S$ and an integer $M$, both only depending on $H_1 \to H$, $S$, and $F$, and a geometric lift $\tilde{\rho}_{\lambda} \colon \gal{F_0, P \cup S_{\lambda}} \to H_1(E_{\lambda})\cdot C^{\vee}[M]$ such that $\tilde{\rho}_{\lambda}$ has good reduction outside $P$. (The assertion that $\tilde{\rho}_{\lambda}$ is valued in the subset $H_1(E_{\lambda}) \cdot C^\vee[M]$ of $\tH(\oE_{\lambda})$ follows from the explicit description of $\tilde{\rho}_{\lambda}$, since $\rho_{\lambda}'$ lands in $H_1(E_{\lambda})$, and $\psi_{\lambda}'$ is trivial.) For all $\lambda$ for which $S_{\lambda} \cap P= \emptyset$, $\nu(\tilde{\rho}_{\lambda}) \colon \gal{F_0, P \cup S_{\lambda}} \to \tZ^\vee[M]$ is also unramified at $S_{\lambda}$, and all such characters are trivialized by a common finite extension $F_1/F_0$. For the finite number of $\lambda$ such that $S_{\lambda} \cap P$ is non-empty, we can again trivialize the possible $\nu(\tilde{\rho}_{\lambda})$ by restricting to a common finite extension $F_2/F_0$. Taking $F'= F_1 F_2$, all $\tilde{\rho}_{\lambda}|_{\gal{F', P \cup S_{\lambda}}}$ land in $H_1(E_{\lambda})$, proving the corollary.
\endproof
\begin{rmk}
For Wintenberger's result, take $E= \Q$. He also shows (\cite[2.3.6]{wintenberger:relevement}) that there is a second finite extension $F''/F'$ (only depending on $H_1 \to H$, $F$, and $S$) such that any two lifts $\rho'_{\lambda}$ as in the corollary become equal after restriction to $\gal{F''}$. This refinement similarly follows in our set-up, but there is no need to repeat Wintenberger's argument.
\end{rmk}
Here is the more general version with multiplicative-type kernels. Note that, as with Theorem \ref{main}, but unlike Corollary \ref{wint}, it makes use of a ``Hodge symmetry" hypothesis.
\begin{cor}\label{multtypecor}
Let $H' \to H$ be a surjection of linear algebraic groups over $E$ whose kernel is central and of multiplicative type. Let $F$ be a number field, and let $S$ be a finite set of places of $F$ containing the archimedean places. Fix a set of co-characters $\{\mu_{\tau}\}_{\tau \colon F \into \oE}$ as in part (2) of Definition \ref{ramcomp}, and moreover assume that each $\mu_{\tau}$ lifts to a co-character of $H'$. 

Then there exist a finite set of places $P \supset S$, and a finite extension $F'/F$, such that any geometric representation $\rho_{\lambda} \colon \gal{F, S \cup S_{\lambda}} \to H(\oE_{\lambda})$ having good reduction outside $S$, and whose Hodge-Tate co-characters arise from the set $\{\mu_{\tau}\}_{\tau \colon F \into \oE}$ via some embedding $\oE \into \oE_{\lambda}$, admits a geometric lift $\tilde{\rho}_{\lambda} \colon \gal{F', P \cup S_{\lambda}} \to H'(\oE_{\lambda})$ having good reduction outside $P$.

In particular, if $\{\rho_{\lambda} \colon \gal{F, S \cup S_{\lambda}} \to H(\oE_{\lambda})\}_{\lambda}$ is a ramification-compatible system with Hodge co-character $\{\mu_{\tau}\}_{\tau \colon F \into \oE}$, then there exist a finite set of places $P \supset S$, a finite extensions $F'/F$, and lifts $\tilde{\rho}_{\lambda} \colon \gal{F', P \cup S_{\lambda}} \to H'(\oE_{\lambda})$ such that $\{\tilde{\rho}_{\lambda}\}_{\lambda}$ is a ramification-compatible system.
\end{cor}
\proof
As in the proofs of Theorem \ref{main} and Corollary \ref{wint}, we construct an isogeny complement $H_1 \subset H$ to $\ker (H' \to H)$, as well as an enlargement $\tH \supset H'$ surjecting onto $H$ with a central torus kernel. We then run the argument of Theorem \ref{main}, starting from lifts $\{\mu_{\tau}'\}$ to $H'$ of the Hodge co-characters: the Hecke character $\psi$ (in the notation of that proof) then constructed has $\lambda$-adic realizations that push-forward to finite-order characters $\psi_{\lambda} \colon \gal{F, T \cup S_{\lambda}} \to \tH/H'(\oE_{\lambda})$, and from here it is easy to proceed: we omit the details, since the argument will by now be familiar. 
\endproof
\bibliographystyle{amsalpha}
\bibliography{biblio.bib}

\end{document}